\documentclass{article}
\PassOptionsToPackage{hyphens}{url}
\usepackage[lang=american]{ems-journal}
\usepackage{mathtools,microtype}
\numberwithin{equation}{section}
\theoremstyle{plain}
\newtheorem{theorem}{Theorem}[section]
\newtheorem{proposition}[theorem]{Proposition}
\newtheorem{lemma}[theorem]{Lemma}

\theoremstyle{definition}

\theoremstyle{remark}
\newtheorem{remark}[theorem]{Remark}
\newcommand{\R}{\mathbb R}
\newcommand{\D}{\mathbb D}
\newcommand{\Sone}{\mathbb S^1}
\newcommand{\Q}{\mathcal Q}

\begin{document}
\title{An Exposition of OpenAI's Planar Hot Spots Theorem}
\hypersetup{pdftitle={An Exposition of OpenAI's Planar Hot Spots Theorem},pdfauthor={Rodrigo Banuelos}}
\titlemark{OpenAI's Planar Hot Spots Theorem}
\date{November 18, 2026}
\emsauthor{1}{\givenname{Rodrigo}\surname{Ba\~nuelos}}{R.~Ba\~nuelos}
\Emsaffil{1}{
  \department{Department of Mathematics}
  \organisation{Purdue University}
  \zip{47907}
  \city{West Lafayette, IN}
   % \zip{47907}
\country{USA}
  \affemail{banuelos@purdue.edu}}
% Add classification and keywords here if needed for submission.
\begin{abstract}We give an exposition of OpenAI's proof that first positive
Neumann eigenfunctions on bounded smooth simply connected planar
domains have no interior critical points. The principal ideas,
techniques, and proof structure are those of the OpenAI preprint.
Our contributions relative to that proof include an alternative
disk Green-kernel regularization argument, a general formulation
of the resolvent step, and a quantitative gradient estimate.
We also interpret the boundary kernels through weighted exit
distributions and boundary trace processes.
\end{abstract}
%\maketitle
\maketitle
\thispagestyle{empty}
%\date{November 18, 2026}
%\begin{center}
%\small November 18, 2026
%\end{center}

\tableofcontents

\section{Introduction and statement of the theorem}
\label{sec:introduction}

Let $G\subset\R^2$ be a bounded simply connected domain with
$C^\infty$ boundary. We write $\Lambda_G^{\mathrm N}$ for the first
positive Neumann eigenvalue of $-\Delta$ and
\[
 E_1(G)=\ker(-\Delta_G-\Lambda_G^{\mathrm N})
\]
for its eigenspace. Here $-\Delta_G$ is the nonnegative self-adjoint
operator associated with the form $v\mapsto\int_G|\nabla v|^2$ on
$H^1(G)$.

The theorem studied in this paper is due to OpenAI \cite{OpenAI}.
It strengthens the strict hot spots conclusion by excluding every
interior critical point, including critical points that are not local
extrema.

\begin{theorem}[OpenAI {\cite{OpenAI}}]
\label{thm:smooth}
Let $G\subset\R^2$ be a bounded simply connected $C^\infty$ domain.
Every nonzero $u\in E_1(G)$ satisfies
\[
 \nabla u(x)\ne0\qquad(x\in G).
\]
In particular,
\begin{equation}\tag{HS1}\label{eq:HS1}
 \min_{\partial G}u<u(x)<\max_{\partial G}u
 \qquad(x\in G).
\end{equation}
No simplicity assumption is imposed on $\Lambda_G^{\mathrm N}$.
\end{theorem}

The strict inequalities follow because $u$ is smooth up to the boundary
and an interior maximum or minimum would be a critical point.
We call \eqref{eq:HS1} the strict hot spots property. The weaker
property (HS2) requires only that the maximum and minimum over the
closure be attained on the boundary. It allows equality at an interior
point as well. These formulations and their early history are
discussed in \cite[Section~1]{BB}.

\paragraph{Relation to the OpenAI \cite{OpenAI} preprint and the companion paper \cite{LipschitzCompanion}.}
Theorem~\ref{thm:smooth} is due to OpenAI \cite{OpenAI}.
Sections~\ref{sec:preliminaries}--\ref{sec:hilbert-finish}
give a detailed exposition of its proof in
\cite[Sections~3--6]{OpenAI}, including the argument for multiple
eigenvalues. The principal arguments, techniques, and their
organization remain those of the OpenAI preprint. References at the
statements identify the results and arguments taken from that work.

The OpenAI proof combines Rohleder's variational principle
\cite{Rohleder2021,Rohleder2025}, a reciprocal-kernel inequality,
and Schoenberg's Hilbert-space realization \cite{Schoenberg}.
The construction of the two positive boundary kernels $K$ and $N$
and the proof that
$K(p,s)K(p,t)/N(s,t)$ is conditionally negative semidefinite
are taken from \cite{OpenAI}. Its boundary multiplier argument
converts the Hilbert-space realization into a sum of nonnegative
energies. A hypothetical interior critical point makes every energy
vanish. The final step uses the bound $\dim E_1(G)\le2$ and boundary
unique continuation, as in the OpenAI proof.

The new contributions developed here, relative to \cite{OpenAI},
concern both the proof and its consequences. We give an alternative
proof of the disk Green-kernel
regularization using hyperbolic distance and Kaluza's sign criterion
\cite{BariczVestiVuorinen2011}. Proposition~\ref{prop:resolvent-preservation}
formulates the positive-quadrature argument of
\cite[Proposition~5.5]{OpenAI} for general continuous kernels,
separating this step from the particular Green kernel of the disk.
The finite-pole
normalization is classical complete Nevanlinna--Pick theory
\cite{AglerMcCarthy2000}. We do not claim that this normalization or
the abstract resolvent formulation is a new result in that theory.

Proposition~\ref{prop:quantitative-gradient} derives a quantitative
gradient estimate from the summed multiplier identity. This is an
additional consequence of the proof. Its convex trace constant is
uniform under fixed inner and outer radius bounds,
while a full uniform gradient estimate remains unresolved here.
These additions do not change the source of the smooth-domain
theorem or the architecture of its proof.

Probabilistic arguments have played a substantial role in the hot spots
literature. Couplings of reflecting Brownian motions underlie results
of Ba\~nuelos and Burdzy \cite{BB}, Pascu \cite{Pascu2002}, and Atar and
Burdzy \cite{AB}. Survival-probability comparisons for Brownian motion
with killing and reflection are central to
\cite{BP2004,BPP}. These successes motivate asking whether stochastic
analysis can explain, or eventually replace, the central kernel
comparison in the smooth-domain proof.

Appendix~\ref{sec:probability} is a preparatory step toward this question.
It identifies a common stochastic structure behind the boundary
kernels. After normalization, $K$ is the exit density
of a gauge-transformed diffusion. The weighted Douglas identity in
Proposition~\ref{prop:weighted-douglas} identifies $N$ as the jump
kernel of the boundary trace of its reflecting version. It expresses
the analytic boundary energy through a Markov process and holds for
the general bounded nonnegative subcritical potentials considered
in the appendix. These identifications use classical Feynman--Kac
and Dirichlet-form theory \cite{ChungZhao1995,FHY2004,FOT2011}.

The probabilistic formulation also isolates a specific obstacle to
a stochastic proof. The central reciprocal kernel has the
Na\"im interpretation given in Remark~\ref{rem:naim}, but its
conditional negative semidefiniteness does not follow from positivity
or the Markov property alone. The higher-dimensional ball calculation
in Remark~\ref{rem:topology-dimension} demonstrates this limitation.
Thus the appendix identifies both the stochastic objects associated
with the proof and the additional geometric inequality that a
probabilistic comparison would need to establish. This gives a
precise target for further stochastic approaches. The appendix does not
provide such a comparison or a probabilistic proof of the theorem.
Its identities supplement the exposition and are not used in the
analytic proof of Theorem~\ref{thm:smooth}.

The notes \cite{Banuelos2026} and \cite{BurdzyNote} showed independently
that the absence of interior local extrema on smooth domains, taken
as a hypothesis, implies a non-strict hot spots conclusion of type
(HS2) for less regular simply connected domains when the first
positive eigenvalue is simple. Both use conformal exhaustion.
The note \cite{Banuelos2026} also treats unbounded simply connected
domains of finite area.

The companion paper \cite{LipschitzCompanion} uses
Theorem~\ref{thm:smooth} to obtain (HS1) for every bounded simply
connected planar Lipschitz domain, without a simplicity assumption.
Its treatment of multiple eigenvalues combines the nodal and
multiplicity results for Lipschitz domains with a thin-strip
construction selecting each prescribed eigendirection
\cite[Lemma~2.3, Corollary~2.4, and Section~3]{LipschitzCompanion}.
Conformal exhaustion then gives smooth approximations
\cite[Lemma~4.3]{LipschitzCompanion}. The Lipschitz extension and these
approximation arguments have no counterpart in \cite{OpenAI}.

The companion's nodal analysis uses two classical ingredients also present in
\cite[Section~2]{OpenAI}. The two-nodal-domain and boundary-sign
arguments are those of \cite[Lemmas~2.3 and~2.4]{OpenAI}, recalled here
in Lemma~\ref{lem:smooth-multiplicity}. The smooth multiplicity bound
is proved in \cite[Proposition~2.5]{OpenAI}. Here we cite its classical
form due to Nadirashvili \cite{Nadirashvili1986} and recalled in
\cite[Proposition~2.5]{BB}. The companion establishes the corresponding
nodal and multiplicity statements for Lipschitz domains.

The strict conclusion in the companion paper follows from its
approximation argument and the exclusion of critical hypersurfaces
by Judge and Mondal
\cite[Proposition~14, p.~160]{JudgeMondalCritical2025}.
The planar argument, using the weak eigenvalue equation, is given in
\cite[Lemma~4.4]{LipschitzCompanion}. The proof of
Theorem~\ref{thm:smooth} in the present paper is independent of the
Lipschitz extension. The appendix's connections with that extension
are explanatory and do not enter the proof.

\paragraph{Historical remarks.}
The hot spots problem goes back to Rauch \cite{Rauch1975}. The
literature is extensive, and we recall only a few results relevant
to the present paper. In addition to the probabilistic results
mentioned above, Jerison and Nadirashvili \cite{JN} proved the strict
hot spots property for convex domains with two orthogonal symmetry
axes. Judge and Mondal \cite{JudgeMondal}, together with the erratum cited there,
proved the result for Euclidean triangles.

Deng, Jiang, and Yang \cite{DengJiangYang2026} prove strict boundary
attainment of the maximum modulus on bounded convex planar domains.
Their November~2, 2026 version states that its proof is independent of
the methods of \cite{OpenAI}. The maximum-modulus conclusion for a
fixed eigenfunction does not by itself imply both inequalities in
\eqref{eq:HS1}. Domains with holes and higher dimensions require
different hypotheses, as the counterexamples in
\cite{BW1999,AtarBurdzy2026,deDiosPont2024} show.
Remark~\ref{rem:topology-dimension} describes the corresponding
obstructions in the present argument.

\paragraph{Organization.}
Section~\ref{sec:preliminaries} gives the spectral and variational
framework. Section~\ref{sec:boundary-reciprocal} constructs the
boundary kernels and proves the multiplier identity.
Section~\ref{sec:reciprocal} proves the reciprocal-kernel inequality.
Section~\ref{sec:hilbert-finish} completes the proof and develops its
quantitative consequences and limitations.
Appendix~\ref{sec:probability} gives the probabilistic interpretations
and concludes with their relation to the Lipschitz approximation.

\section{Variational framework on smooth domains}
\label{sec:preliminaries}

Throughout the proof, $G\subset\R^2$ is bounded and simply connected
with $C^\infty$ boundary. All functions and Sobolev spaces are real.
A smooth domain means that, near every boundary point and after a
rigid motion, the domain lies on one side of the graph of a
$C^\infty$ function. The same definition with a Lipschitz graph is
used for the auxiliary domains in the appendix.

\subsection{Notation and spectral facts}
\label{subsec:conventions}

The Neumann eigenvalues are indexed from zero and counted with
multiplicity. For any bounded Lipschitz domain $\Omega$ we write
\[
 0=\Lambda_{\Omega,0}^{\mathrm N}<\Lambda_{\Omega,1}^{\mathrm N}
 \le\Lambda_{\Omega,2}^{\mathrm N}\le\cdots,
 \qquad \Lambda_\Omega^{\mathrm N}=\Lambda_{\Omega,1}^{\mathrm N},
\]
and denote the lowest Dirichlet eigenvalue by
$\lambda_\Omega^{\mathrm D}$. The definition of $E_1(G)$ extends to
$E_1(\Omega)$ in the same way. A first eigenfunction is a nonzero
element of this eigenspace. A nodal domain is a connected component
of its nonzero set. The variational characterizations are
\begin{align}
 \Lambda_\Omega^{\mathrm N}
 &=\inf_{\substack{v\in H^1(\Omega)\setminus\{0\}\\\int_\Omega v=0}}
 \frac{\int_\Omega|\nabla v|^2}{\int_\Omega v^2},
 \label{eq:neumann-minmax}\\
 \lambda_\Omega^{\mathrm D}
 &=\inf_{v\in H^1_0(\Omega)\setminus\{0\}}
 \frac{\int_\Omega|\nabla v|^2}{\int_\Omega v^2}.
 \label{eq:dirichlet-minmax}
\end{align}
We use the classical strict comparison
\begin{equation}\label{eq:strict-spectral-gap}
 \Lambda_\Omega^{\mathrm N}<\lambda_\Omega^{\mathrm D},
\end{equation}
which is a special case of Filonov's theorem
\cite[p.~414]{Filonov}. Neither simple connectedness nor simplicity
of the Neumann eigenvalue is required. For the earlier planar results
of P\'olya \cite{Polya1952} and Szeg\H{o} \cite{Szego1954}, see
\cite[p.~413]{Filonov}.

We write $\D=\{z\in\mathbb C:|z|<1\}$ and $\Sone=\partial\D$.
Boundary integrals use arclength, so that $\int_{\Sone}ds=2\pi$.
The symbol $\Phi$ always denotes the conformal map onto $G$.
Conformal maps associated with the companion paper's approximation
are denoted by $F$ and defined locally in the appendix.

\subsection{Regularity and the conformal map}
\label{subsec:regularity}

Standard elliptic boundary
regularity \cite[Chapter~III, Theorems~6.1 and~6.5]{Showalter} gives
\begin{equation}\label{eq:eigenfunction-regularity}
 v\in E_1(G)\quad\Longrightarrow\quad
 v\in C^\infty(\overline G),\qquad
 -\Delta v=\Lambda_G^{\mathrm N}v\ \hbox{in }G,\qquad
 \partial_\nu v=0\ \hbox{on }\partial G.
\end{equation}
Here $C^\infty(\overline G)$ means smoothness up to the boundary.

Let
\begin{equation}\label{eq:conformal-map}
 \Phi:\D\longrightarrow G
\end{equation}
be a conformal bijection. By the Kellogg--Warschawski theorem
\cite[Theorems~3.5--3.6]{Pommerenke}, $\Phi$ extends to a
$C^\infty$ diffeomorphism of the closures, and
\begin{equation}\label{eq:conformal-nondegeneracy}
 \min_{\overline\D}|\Phi'|>0.
\end{equation}
Thus $\Phi'$ and $1/\Phi'$ are smooth and bounded on
$\overline\D$.

\subsection{Nodal domains and multiplicity}
\label{subsec:smooth-multiplicity}

The multiplicity bound below is the classical smooth-boundary result
of Nadirashvili \cite{Nadirashvili1986}. See also
\cite[Proposition~2.5]{BB}. We use this known bound as a spectral
input. The nodal-domain and boundary-sign arguments are recalled for
the probabilistic interpretation and follow
\cite[Lemmas~2.3 and~2.4]{OpenAI}.

\begin{lemma}[Smooth-domain nodal facts and multiplicity]
\label{lem:smooth-multiplicity}
Every nonzero $u\in E_1(G)$ has exactly two nodal domains and takes
both signs on $\partial G$. Moreover,
\[
 \dim E_1(G)\le2.
\]
\end{lemma}

\begin{proof}
Since $\int_Gu=0$, there are at least two nodal domains. If there
were at least three, choose a positive domain $G_1$ and a negative
domain $G_2$. The nodal restrictions $u_i=u\mathbf1_{G_i}$ belong to
$H^1(G)$ and satisfy $\nabla u_i=\mathbf1_{G_i}\nabla u$.
Indeed, for $c>0$, the function $(u-c)^+\mathbf1_{G_1}$ vanishes
near every interior point of $\partial G_1$. Its weak gradient is
$\mathbf1_{G_1\cap\{u>c\}}\nabla u\in L^2(G)$. Dominated convergence
gives convergence in $H^1(G)$ to $u\mathbf1_{G_1}$ as $c\downarrow0$.
Apply the same argument to $-u$ on $G_2$.
Testing the weak equation with $u_i$ gives
\[
 \int_G|\nabla u_i|^2=\Lambda_G^{\mathrm N}\int_Gu_i^2.
\]
The nonzero function
\[
 v=\left(\int_Gu_2\right)u_1
       -\left(\int_Gu_1\right)u_2
\]
has mean zero and Rayleigh quotient $\Lambda_G^{\mathrm N}$.
Thus $v\in E_1(G)$, but $v$ vanishes on the third nodal domain,
contrary to interior unique continuation.

If $u\ge0$ on $\partial G$, then the nonzero negative part $u^-$
belongs to $H^1_0(G)$. Its Rayleigh quotient is
$\Lambda_G^{\mathrm N}$, contradicting
\eqref{eq:strict-spectral-gap}. Applying the argument to $-u$ gives
both boundary signs. The final assertion is the cited multiplicity
theorem of Nadirashvili for smooth simply connected planar domains.
\end{proof}

\subsection{Rohleder's variational principle}
\label{sec:vector-energy}

Let $\nu$ denote the outward unit normal, and put
\begin{align}
 \mathcal X(G)&=\{X\in H^1(G;\R^2):
                    \operatorname{Tr}X\cdot\nu=0\},
                       \label{eq:tangent-vector-space}\\
 \operatorname{curl}X&=\partial_1X_2-\partial_2X_1,
 \qquad
 \Q(X)=\int_G\big((\operatorname{div}X)^2+
                         (\operatorname{curl}X)^2\big)\,dx.
                       \label{eq:vector-energy}
\end{align}
The following statement is Rohleder's variational principle. In
\cite[Theorem~1.1]{Rohleder2021} it is stated in the equivalent
curvature form given by Lemma~\ref{lem:curvature} below, for piecewise
smooth domains with convex corners and with the opposite sign of
$\kappa$. It also follows from the spectral description and variational
characterization of the divergence--curl form in
\cite[Proposition~3.2 and Theorem~3.4]{Rohleder2025}. The
formulation and the proof below, including the equality case, follow
\cite[Proposition~3.1]{OpenAI}.

\begin{proposition}[Variational principle of Rohleder {\cite[Theorem~1.1]{Rohleder2021}}]
\label{prop:vector}
For every $X\in\mathcal X(G)$,
\begin{equation}\label{eq:vector-poincare}
 \Q(X)\ge\Lambda_G^{\mathrm N}\int_G|X|^2\,dx,
\end{equation}
and
\begin{equation}\label{eq:vector-equality}
 \Q(X)=\Lambda_G^{\mathrm N}\int_G|X|^2\,dx
 \quad\Longleftrightarrow\quad
 X=\nabla v\quad\hbox{for some }v\in E_1(G).
\end{equation}
\end{proposition}

No assumption on the multiplicity of $\Lambda_G^{\mathrm N}$ is made.

\begin{proof}
Write $\nabla^\perp=(-\partial_2,\partial_1)$ and
$\|\cdot\|_2=\|\cdot\|_{L^2(G)}$.

\emph{Decomposition.}
Let $v$ and $\phi$ be the solutions of
\[
 \Delta v=\operatorname{div}X\ \text{in }G,\quad
 \partial_\nu v=0\ \text{on }\partial G,\quad\int_Gv\,dx=0,
 \qquad\qquad
 \Delta\phi=\operatorname{curl}X\ \text{in }G,\quad
 \phi=0\ \text{on }\partial G.
\]
The Neumann problem is solvable because
$\int_G\operatorname{div}X\,dx=\int_{\partial G}\operatorname{Tr}X\cdot\nu\,ds=0$.
Both right-hand sides belong to $L^2(G)$, so $v,\phi\in H^2(G)$ by
elliptic boundary regularity
\cite[Chapter~III, Theorems~6.1 and~6.4]{Showalter}. The field
$\nabla^\perp\phi$ is divergence-free, satisfies
$\operatorname{curl}\nabla^\perp\phi=\Delta\phi$, and is tangent to
$\partial G$ because $\phi=0$ on $\partial G$. Hence
\[
 X_0=X-\nabla v-\nabla^\perp\phi\in H^1(G;\R^2)
\]
satisfies
\[
 \operatorname{div}X_0=0,\qquad\operatorname{curl}X_0=0,\qquad
 \operatorname{Tr}X_0\cdot\nu=0.
\]
We show that $X_0=0$. Write $X_0=(X_{0,1},X_{0,2})$. The first two
equations are the Cauchy--Riemann equations for $X_{0,1}-iX_{0,2}$, and
they hold almost everywhere in $G$. By Weyl's lemma,
$X_{0,1}-iX_{0,2}$ is holomorphic in $G$. Since $G$ is simply connected,
this function has a holomorphic primitive, and the real part $\eta$ of
the primitive satisfies $\nabla\eta=X_0$ in $G$. The function
$\eta\circ\Phi$ is smooth in $\D$, and by conformal invariance its
Dirichlet integral equals $\|X_0\|_2^2<\infty$. For
$\tfrac12<\varrho\le1$ and $f\in H^1(\varrho\D)$ we have
\[
 \|f\|_{L^2(\varrho\D)}^2\le C\bigl(\|\nabla f\|_{L^2(\varrho\D)}^2
 +\|f\|_{L^2(\frac12\D)}^2\bigr)
\]
with a constant $C$ independent of $\varrho$: the case $\varrho=1$ is a
standard form of the Poincar\'e inequality on $\D$, and the general
case follows by applying it to $x\mapsto f(\varrho x)$. We apply this
inequality to $f=\eta\circ\Phi$, which is smooth on
$\varrho\overline\D$ for $\varrho<1$. Letting $\varrho\uparrow1$,
monotone convergence gives $\eta\circ\Phi\in L^2(\D)$, and therefore
$\eta\circ\Phi\in H^1(\D)$. Since $\Phi$ is a diffeomorphism of the
closures, $\eta\in H^1(G)$. Green's formula now gives
\[
 \|X_0\|_2^2=\int_GX_0\cdot\nabla\eta\,dx
 =-\int_G(\operatorname{div}X_0)\,\eta\,dx
  +\int_{\partial G}(\operatorname{Tr}X_0\cdot\nu)\,\eta\,ds=0.
\]
Thus
\begin{equation}\label{eq:vector-helmholtz}
 X=\nabla v+\nabla^\perp\phi.
\end{equation}

\emph{The inequality.}
The two summands in \eqref{eq:vector-helmholtz} are orthogonal in
$L^2(G;\R^2)$, because $\nabla^\perp\phi$ is divergence-free and
tangent:
$\int_G\nabla v\cdot\nabla^\perp\phi\,dx
=-\int_Gv\operatorname{div}\nabla^\perp\phi\,dx
+\int_{\partial G}v\,\nabla^\perp\phi\cdot\nu\,ds=0$. Hence
\[
 \|X\|_2^2=\|\nabla v\|_2^2+\|\nabla\phi\|_2^2,\qquad
 \Q(X)=\|\Delta v\|_2^2+\|\Delta\phi\|_2^2.
\]
Integration by parts and the variational characterizations
\eqref{eq:neumann-minmax} and \eqref{eq:dirichlet-minmax} give
\begin{equation}\label{eq:vector-two-rayleigh}
 \|\nabla v\|_2^2=-\int_Gv\,\Delta v\,dx
 \le\|v\|_2\|\Delta v\|_2
 \le(\Lambda_G^{\mathrm N})^{-1/2}\|\nabla v\|_2\|\Delta v\|_2,
\end{equation}
and in the same way
$\|\nabla\phi\|_2^2\le(\lambda_G^{\mathrm D})^{-1/2}\|\nabla\phi\|_2\|\Delta\phi\|_2$.
Therefore
\begin{equation}\label{eq:vector-gap-remainder}
 \Q(X)\ge\Lambda_G^{\mathrm N}\|\nabla v\|_2^2
        +\lambda_G^{\mathrm D}\|\nabla\phi\|_2^2
 =\Lambda_G^{\mathrm N}\|X\|_2^2
  +(\lambda_G^{\mathrm D}-\Lambda_G^{\mathrm N})\|\nabla\phi\|_2^2 .
\end{equation}
Since $\Lambda_G^{\mathrm N}<\lambda_G^{\mathrm D}$ by
\eqref{eq:strict-spectral-gap}, this proves
\eqref{eq:vector-poincare}.

\emph{The equality case.}
Suppose that equality holds in \eqref{eq:vector-poincare}. Then
\eqref{eq:vector-gap-remainder} gives $\nabla\phi=0$, hence $\phi=0$
and $X=\nabla v$, and also
$\|\Delta v\|_2^2=\Lambda_G^{\mathrm N}\|\nabla v\|_2^2$. If
$\nabla v=0$, then $v=0\in E_1(G)$. Otherwise both inequalities in
\eqref{eq:vector-two-rayleigh} are equalities. Equality in the
Cauchy--Schwarz step gives $-\Delta v=cv$ with a constant $c>0$, and
equality in the second step gives $c=\Lambda_G^{\mathrm N}$. Since
$\partial_\nu v=0$, we get $v\in E_1(G)$. Conversely, if
$X=\nabla\widetilde v$ with $\widetilde v\in E_1(G)$, then
$X\in\mathcal X(G)$ by \eqref{eq:eigenfunction-regularity}, and
\[
 \Q(X)=\|\Delta\widetilde v\|_2^2
 =(\Lambda_G^{\mathrm N})^2\|\widetilde v\|_2^2
 =\Lambda_G^{\mathrm N}\|X\|_2^2.
\]
\end{proof}

Let $\tau$ be the counterclockwise unit tangent and define
\begin{equation}\label{eq:curvature-convention}
 \kappa=\det(\tau,\partial_\tau\tau).
\end{equation}
Our convention gives $\kappa>0$ on a circle and is opposite to that
in \cite{Rohleder2021}. The following identity gives the curvature
form of the variational principle. See \cite{Rohleder2024} and
\cite[Remark~3.6, p.~5563]{Rohleder2025}.

\begin{lemma}[Curvature form of the energy; {\cite[Lemma~3.2]{OpenAI}}]\label{lem:curvature}
For $X\in\mathcal X(G)$,
\begin{equation}\label{eq:curvature-energy}
 \Q(X)=\int_G|\nabla X|^2\,dx
          +\int_{\partial G}\kappa|\operatorname{Tr}X|^2\,ds,
 \qquad |\nabla X|^2=\sum_{i,j=1}^2|\partial_iX_j|^2.
\end{equation}
\end{lemma}

\begin{proof}
For a smooth tangent field, write $X=a\tau$ on $\partial G$.
Expansion and Stokes' formula give
\[
 \Q(X)-\int_G|\nabla X|^2\,dx
 =2\int_G\det\nabla X\,dx
 =\int_{\partial G}\det(X,\partial_\tau X)\,ds
 =\int_{\partial G}\kappa a^2\,ds.
\]
For general $X\in\mathcal X(G)$, write $\operatorname{Tr}X=b\tau$.
Let $\mathcal L$ be the componentwise harmonic Dirichlet lifting.
It is bounded from $H^{1/2}(\partial G)$ to $H^1(G)$ and takes
smooth traces to functions in $C^\infty(\overline G)$.
Choose smooth $b_j\to b$ in $H^{1/2}(\partial G)$ and
$Z_j\in C_c^\infty(G;\R^2)$ converging to
$X-\mathcal L(b\tau)\in H^1_0(G;\R^2)$. Then
\[
 X_j=\mathcal L(b_j\tau)+Z_j\longrightarrow X\quad\text{in }H^1(G;\R^2),
 \qquad \operatorname{Tr}X_j=b_j\tau.
\]
Applying the smooth identity to $X_j$ and passing to the limit
proves the assertion.
\end{proof}

\section{Boundary kernels and the multiplier identity}
\label{sec:boundary-reciprocal}

\label{sec:kernels}

We construct the kernels $K$, $R$, and $N$, following
\cite[Section~4]{OpenAI} and giving the details. The construction uses
only that $a$ is bounded and nonnegative and satisfies
\eqref{eq:subcritical}. It does not use that
$a=\Lambda_G^{\mathrm N}|\Phi'|^2$ (see
\cite[Sections~1 and~4.1]{OpenAI}). We need this generality in
Section~\ref{sec:reciprocal}, where $a$ is replaced by
$a\mathbf1_{\{|x|\le r_n\}}$. Section~\ref{sec:probability} expresses
these objects through Brownian motion.

\subsection{The disk resolvent}

Let $dq=a\,dx$, where $a\in L^\infty(\D)$ is nonnegative, and suppose that
\begin{equation}
 \int_\D |\psi|^2\,dq
 \le d\int_\D |\nabla\psi|^2\,dx,
 \qquad \psi\in H^1_0(\D),\qquad 0<d<1.
 \label{eq:subcritical}
\end{equation}
The symbols $\langle\cdot,\cdot\rangle_q$ and $\|\cdot\|_q$ denote the
inner product and norm in $L^2(q)$. The Dirichlet Green kernel and
Poisson kernel of the disk are
\begin{align}
 \mathcal G(x,y)
 &=\frac1{2\pi}\log\left|\frac{1-x\overline y}{x-y}\right|
 =\frac1{4\pi}\log\left(1+
       \frac{(1-|x|^2)(1-|y|^2)}{|x-y|^2}\right),
       \label{eq:disk-green}\\
 P_s(x)&=\frac{1-|x|^2}{2\pi|s-x|^2},
 \qquad x\in\D,\quad s\in\Sone.
 \label{eq:disk-poisson}
\end{align}
Thus $\mathcal G(x,y)>0$ for $x\ne y$, and
$\int_{\Sone}P_s(x)\,ds=1$.  We set
$\mathcal G(x,x)=+\infty$. The value on the diagonal has no effect on
integration against $q$. For $x,y\in\D$ with $x\ne y$ and $s\in\Sone$,
we shall use the bounds
\begin{equation}\label{eq:green-poisson-bounds}
 P_s(y)\le\frac1{\pi|y-s|},\qquad
 \mathcal G(x,y)\le\frac1{2\pi}\log\frac2{|x-y|}.
\end{equation}

\begin{lemma}[The disk resolvent; cf.\ {\cite[Lemma~4.1 and~(4.5)]{OpenAI}}]
\label{lem:resolvent}
Let $a$, $q$, and $d$ be as in \eqref{eq:subcritical}.
\begin{enumerate}
\item[\textup{(i)}] The operator
\begin{equation}\label{eq:green-operator}
 Tf(x)=\int_\D\mathcal G(x,y)f(y)\,dq(y)
\end{equation}
is self-adjoint, nonnegative and positivity preserving on $L^2(q)$,
with $\|T\|\le d$.
\item[\textup{(ii)}] Its resolvent $S=(I-T)^{-1}$
is also self-adjoint, nonnegative and positivity preserving, and
\begin{equation}\label{eq:L2-resolvent-bound}
 \|S\|_{L^2(q)\to L^2(q)}\le\frac1{1-d}.
\end{equation}
\item[\textup{(iii)}] We have
\begin{equation}\label{eq:green-evaluation-bound}
 C_a:=\sup_{x\in\D}\|\mathcal G(x,\cdot)\|_q<\infty,
 \qquad \|Tf\|_\infty\le C_a\|f\|_q\quad(f\in L^2(q)).
\end{equation}
\item[\textup{(iv)}] For Lipschitz boundary data $h:\Sone\to\R^m$, let
\[
 Ah(x)=\int_{\Sone}P_s(x)h(s)\,ds
\]
be its Poisson integral. There is a unique
$W_h\in H^1(\D;\R^m)$ that satisfies $-\Delta W_h=aW_h$
weakly and has trace $h$. It is continuous on $\overline\D$, and
\begin{equation}\label{eq:solution-resolvent}
 W_h=Ah+TSAh,\qquad W_h=SAh\quad\text{in }L^2(q).
\end{equation}
\item[\textup{(v)}] Let $\rho$ be a measure on $\D$ with
$0\le\rho\le q$, and let $a_\rho$ be its density. Then
\eqref{eq:subcritical} holds for $\rho$ with the same constant $d$, and
\textup{(i)}--\textup{(iv)} hold with $q$ replaced by $\rho$ and $a$
replaced by $a_\rho$.
\end{enumerate}
\end{lemma}

Parts \textup{(i)} and \textup{(v)} contain \cite[Lemma~4.1]{OpenAI},
and part \textup{(iv)} is \cite[(4.5)]{OpenAI}. We also record the
standard positivity-preserving properties of $T$ and $S$.

\begin{proof}
For $f\in L^2(q)$, $af\in L^2(\D)$ and
$\psi=Tf\in H^2(\D)\cap H^1_0(\D)$ solves $-\Delta\psi=af$.
Hence
\[
 \|Tf\|_q^2\le d\int_\D|\nabla Tf|^2
 =d\langle Tf,f\rangle_q\le d\|Tf\|_q\|f\|_q,
 \qquad
 \langle f,Tk\rangle_q=\int_\D\nabla Tf\cdot\nabla Tk.
\]
Thus $T$ is self-adjoint, nonnegative, and $\|T\|\le d$.
It preserves positivity because $\mathcal G>0$.
The spectral theorem gives $S$ and \eqref{eq:L2-resolvent-bound}.
If $f\ge0$ and $v=TSf$, testing $(-\Delta-a)v=af$ by $-v^-$ gives
\[
 (1-d)\int_\D|\nabla v^-|^2
 \le-\int_\D afv^-\le0.
\]
Hence $Sf=f+v\ge0$. This proves \textup{(i)} and \textup{(ii)}. The
logarithmic bound in \eqref{eq:green-poisson-bounds} and the
boundedness of $a$ give \textup{(iii)}, by the Cauchy--Schwarz
inequality.

For Lipschitz $h$, set $c=SAh$ and $W_h=Ah+Tc$.
Then $W_h=c$ in $L^2(q)$, $-\Delta W_h=aW_h$, and
$\operatorname{Tr}W_h=h$. Since $Tc\in H^2\cap H^1_0$,
$W_h\in H^1(\D)\cap C(\overline\D)$.
The form bound \eqref{eq:subcritical} gives uniqueness by testing
the homogeneous zero-trace equation with its solution. The same
proof applies componentwise to vector data. This proves \textup{(iv)}.
If $0\le\rho\le q$, the density of $\rho$ is bounded by $a$
and its form bound has the same $d$, which gives \textup{(v)}.
\end{proof}

\subsection{The kernels \texorpdfstring{$K$, $R$, and $N$}{K, R, and N}}

In general $P_s\notin L^2(q)$, so we first apply $T$ to it. By
\eqref{eq:green-poisson-bounds},
\[
 \sup_{s\in\Sone}\|P_s\|_{L^{3/2}(\D)}<\infty,\qquad
 \sup_{x\in\D}\|\mathcal G(x,\cdot)\|_{L^3(\D)}<\infty.
\]
H\"older's inequality therefore yields
\begin{equation}\label{eq:poisson-smoothed}
 F_s:=TP_s\in L^\infty(\D)\cap L^2(q),\qquad
 \sup_{s\in\Sone}\|F_s\|_\infty\le C\|a\|_\infty.
\end{equation}
We use the bounded representative
\[
 U_s:=F_s+TSF_s
\]
of $SF_s$. For $p\in\D$ and distinct $s,t\in\Sone$, define
\begin{align}
 K(p,s)&=P_s(p)+U_s(p),\label{eq:kernel-K}\\
 R(s,t)&=\int_\D P_s(x)K(x,t)\,dq(x),\label{eq:kernel-R}\\
 N(s,t)&=N_0(s,t)+R(s,t),\qquad
 N_0(s,t)=\frac1{\pi|s-t|^2}.\label{eq:kernel-N}
\end{align}

\begin{lemma}[The kernels \texorpdfstring{$K$, $R$, and $N$}{K, R, and N}; {\cite[Lemma~4.2]{OpenAI}}]\label{lem:kernels}
The kernels in \eqref{eq:kernel-K}--\eqref{eq:kernel-N} are finite at the
specified points. For fixed $p\in\D$, $K(p,\cdot)$ is positive and
bounded. The kernel $R$ is symmetric, nonnegative and integrable, with
\begin{equation}\label{eq:R-integrable}
 \int_{\Sone}\int_{\Sone}R(s,t)\,ds\,dt
 \le\frac{q(\D)}{1-d}.
\end{equation}
In particular, $N(s,t)\ge N_0(s,t)>0$ for $s\ne t$.
For Lipschitz scalar or vector data $h,k$,
\begin{align}
 W_h(p)&=\int_{\Sone}K(p,s)h(s)\,ds,
 \label{eq:K-representation}\\
 \langle Ah,SAk\rangle_q
 &=\int_{\Sone}\int_{\Sone}R(s,t)h(s)\cdot k(t)\,ds\,dt.
 \label{eq:R-representation}
\end{align}
Moreover, for $p\in\D$ and distinct $s,t\in\Sone$, the kernels have the
series representations
\begin{equation}\label{eq:kernel-chains}
 \begin{aligned}
 K(p,s)&=P_s(p)+\sum_{n=1}^{\infty}(T^nP_s)(p),\\
 R(s,t)&=\sum_{n=0}^{\infty}
              \int_\D P_s(x)(T^nP_t)(x)\,dq(x),
 \end{aligned}
\end{equation}
in which every term is nonnegative. Here, for $n\ge1$, $(T^nP_s)(p)$
denotes the iterated integral
\[
 \int_{\D^n}\mathcal G(p,y_1)\mathcal G(y_1,y_2)\cdots
 \mathcal G(y_{n-1},y_n)P_s(y_n)\,dq(y_1)\cdots dq(y_n),
\]
which equals $(T^{n-1}F_s)(p)$ by Tonelli's theorem.
\end{lemma}

\begin{proof}
Positivity of $T$ and $S$, together with
\eqref{eq:green-evaluation-bound}, gives
\[
 0\le U_s,\qquad
 \|U_s\|_\infty
 \le\|F_s\|_\infty+\frac{C_a}{1-d}\|F_s\|_q.
\]
The bound is uniform in $s$. Since $P_s(p)$ is positive and bounded
as a function of $s$ for fixed $p$, the assertion about $K$ follows.

For $s\ne t$, $\int P_sP_t\,dq<\infty$. Indeed, near either
boundary point one Poisson kernel is bounded and the other is integrable.
Also $\sup_s\int P_s\,dq<\infty$. Hence $R(s,t)<\infty$, and
Tonelli's theorem gives
\begin{equation}\label{eq:R-resolvent-identity}
 R(s,t)=\int_\D P_sP_t\,dq+\int_\D P_sF_t\,dq
                         +\langle F_s,SF_t\rangle_q.
\end{equation}
The first term is symmetric, the second is symmetric by symmetry of
$\mathcal G$, and the third by self-adjointness of $S$.

Rotations act continuously on $L^{3/2}(\D)$, so $s\mapsto P_s$ is
continuous in that space. The estimate used in
\eqref{eq:poisson-smoothed} gives
\[
 \|F_s-F_t\|_\infty
 \le\|a\|_\infty
     \sup_{x\in\D}\|\mathcal G(x,\cdot)\|_{L^3(\D)}
     \|P_s-P_t\|_{L^{3/2}(\D)}\longrightarrow0
 \qquad(t\to s).
\]
Thus $s\mapsto F_s$ is continuous in $L^2(q)$, and
$\int_{\Sone}F_sh(s)\,ds$ is a Bochner integral in that space.
The bounded operator $TS:L^2(q)\to C(\overline\D)$ may be applied
under the integral sign. Its boundedness also gives joint continuity
of $(p,s)\mapsto(TSF_s)(p)$. The other terms in
\eqref{eq:R-resolvent-identity} are measurable by integration, so
$K$ and $R$ are jointly measurable.
Since $\int F_sh\,ds=TAh$ and
$S=I+ST=I+TS$,
\[
 \int_{\Sone}K(p,s)h(s)\,ds
 =Ah(p)+TAh(p)+(TSTAh)(p)
 =Ah(p)+(TSAh)(p)=W_h(p).
\]
For nonnegative scalar $h,k$, Tonelli now yields
\[
 \iint R(s,t)h(s)k(t)\,ds\,dt
 =\int_\D Ah(x)W_k(x)\,dq(x)
 =\langle Ah,SAk\rangle_q.
\]
In particular, taking $h=k=1$ proves directly that
\[
 \iint R(s,t)\,ds\,dt
 =\langle1,S1\rangle_q
 \le\|S\|\|1\|_q^2\le\frac{q(\D)}{1-d}.
\]
For signed scalar data and for vector data,
\eqref{eq:R-representation} follows by linearity. The integrals
converge absolutely because $R$ is integrable and $h$, $k$ are
bounded.
The diagonal values of $R$ may be chosen arbitrarily.

Finally, the Neumann series $S=\sum_{n\ge0}T^n$ converges in operator
norm on $L^2(q)$, because $\|T\|\le d<1$. Since $F_s\in L^2(q)$ and
$\|\mathcal G(p,\cdot)\|_q\le C_a$,
\[
 (TSF_s)(p)=\Bigl\langle\mathcal G(p,\cdot),\sum_{n\ge0}T^nF_s\Bigr\rangle_q
 =\sum_{n\ge0}(T^{n+1}F_s)(p)
\]
for every $p\in\D$, which is the first line of
\eqref{eq:kernel-chains}. For the second line, the first two terms of
\eqref{eq:R-resolvent-identity} are the terms $n=0$ and $n=1$. In the
third term we insert the Neumann series and use
$\langle F_s,T^mF_t\rangle_q=\int_\D P_s\,T^{m+1}F_t\,dq$ for $m\ge0$,
which follows from Tonelli's theorem and the symmetry of $\mathcal G$;
this is the term $n=m+2$. These expansions are valid although $P_s$
may fail to belong to $L^2(q)$, because $F_s=TP_s$ belongs to
$L^2(q)$.
\end{proof}

\subsection{The vector boundary energy}

We now return to the smooth simply connected domain $G$, its first
positive Neumann eigenvalue $\Lambda_G^{\mathrm N}$, and a conformal diffeomorphism
$\Phi:\overline\D\to\overline G$. In the rest of
Section~\ref{sec:kernels},
\begin{equation}\label{eq:conformal-potential}
 dq(x)=\Lambda_G^{\mathrm N}|\Phi'(x)|^2\,dx,
\end{equation}
so that $a=\Lambda_G^{\mathrm N}|\Phi'|^2$,
and $T$, $S$, $W_h$, $K$, $R$, and $N$ are the objects of
Lemmas~\ref{lem:resolvent} and~\ref{lem:kernels} for this measure.
Conformal invariance and \eqref{eq:strict-spectral-gap}, with the notation
$\lambda_G^{\mathrm D}$ of \eqref{eq:dirichlet-minmax}, give
\eqref{eq:subcritical} with $d=\Lambda_G^{\mathrm N}/\lambda_G^{\mathrm D}<1$.  Indeed, for
$\psi\in H^1_0(\D)$ and $\widetilde\psi=\psi\circ\Phi^{-1}$,
\[
 \int_\D|\psi|^2\,dq
 =\Lambda_G^{\mathrm N}\int_G|\widetilde\psi|^2
 \le\frac{\Lambda_G^{\mathrm N}}{\lambda_G^{\mathrm D}}\int_G|\nabla\widetilde\psi|^2
 =\frac{\Lambda_G^{\mathrm N}}{\lambda_G^{\mathrm D}}\int_\D|\nabla\psi|^2.
\]

We call a map $h:\Sone\to\R^2$ \emph{tangent} if
$h(s)\cdot\nu(\Phi(s))=0$ for every $s\in\Sone$. For tangent Lipschitz
maps $h,k$, define
\begin{equation}\label{eq:energy-definition}
 X_h=W_h\circ\Phi^{-1},\qquad
 E(h,k)=\Q(X_h,X_k)-\Lambda_G^{\mathrm N}\int_G X_h\cdot X_k\,dx,\qquad E(h)=E(h,h),
\end{equation}
where $\Q(\cdot,\cdot)$ is the polarization of the divergence--curl form
in Proposition~\ref{prop:vector}. The field $X_h$ is obtained by
composing each component of $W_h$ with $\Phi^{-1}$. It is not the
push-forward of $W_h$ as a vector field.

\begin{lemma}[Nonnegativity and null space of the boundary energy; {\cite[(4.13)]{OpenAI}}]
\label{lem:energy-nullspace}
For every tangent Lipschitz map $h:\Sone\to\R^2$, the field $X_h$
belongs to $\mathcal X(G)$, and
\begin{equation}
 E(h)\ge0,\qquad
 E(h)=0\ \Longleftrightarrow\quad
 X_h=\nabla v\text{ for some }v\in E_1(G).
 \label{eq:energy-nullspace}
\end{equation}
\end{lemma}

\begin{proof}
Since $W_h\in H^1(\D;\R^2)$ and $\Phi$ is a diffeomorphism of the
closures, $X_h\in H^1(G;\R^2)$. By
Lemma~\ref{lem:resolvent}\,\textup{(iv)}, $X_h$ is continuous on
$\overline G$ and $X_h(\Phi(s))=h(s)$, which is orthogonal to
$\nu(\Phi(s))$. Thus $X_h\in\mathcal X(G)$, and
\eqref{eq:energy-nullspace} is Proposition~\ref{prop:vector}.
\end{proof}

\begin{lemma}[Boundary formula for the energy; {\cite[Lemma~4.3]{OpenAI}}]\label{lem:boundary-energy}
Put
\[
 r(s)=\int_{\Sone}R(s,t)\,dt .
\]
Then $r\in L^1(\Sone)$, and for tangent Lipschitz maps
$h,k:\Sone\to\R^2$,
\begin{align}
 E(h,k)={}&\frac12\iint_{\Sone\times\Sone}N(s,t)
       (h(s)-h(t))\cdot(k(s)-k(t))\,ds\,dt\notag\\
 &+\int_{\Sone}
       \bigl(|\Phi'(s)|\kappa(\Phi(s))-r(s)\bigr)
                            h(s)\cdot k(s)\,ds.
 \label{eq:boundary-energy}
\end{align}
Here $\kappa$ has the sign convention of Lemma~\ref{lem:curvature}.
\end{lemma}

\begin{proof}
The function $r$ is integrable by \eqref{eq:R-integrable}.
Put
\[
 \mathcal B(h,k)
   =\int_\D\nabla W_h:\nabla W_k\,dx
                       -\langle W_h,W_k\rangle_q.
\]
Testing $-\Delta W_k=aW_k$ with $W_h-Ah\in H^1_0(\D)$ gives
$\int_\D\nabla W_k:\nabla(W_h-Ah)\,dx=\langle W_k,W_h-Ah\rangle_q$.
Since $Ah$ is harmonic and $W_k-Ak\in H^1_0(\D)$,
$\int_\D\nabla Ah:\nabla(W_k-Ak)\,dx=0$. Together with $W_k=SAk$ in
$L^2(q)$, these identities give
\[
 \mathcal B(h,k)
  =\int_\D\nabla Ah:\nabla Ak\,dx-\langle Ah,SAk\rangle_q.
\]
The polarized Douglas formula is
\begin{equation}\label{eq:harmonic-boundary-energy}
 \int_\D\nabla Ah:\nabla Ak\,dx
 =\frac1{2\pi}\iint_{\Sone\times\Sone}
       \frac{(h(s)-h(t))\cdot(k(s)-k(t))}{|s-t|^2}\,ds\,dt.
\end{equation}
See \cite{Douglas1931} and \cite[(1.1)]{BGPR2023}. The latter formula
is written in the angular variables $\eta,\xi$ of $s=e^{i\eta}$,
$t=e^{i\xi}$, with constant $1/(8\pi)$ and kernel
$\sin^{-2}((\eta-\xi)/2)$; it agrees with
\eqref{eq:harmonic-boundary-energy} because
$|s-t|^2=4\sin^2((\eta-\xi)/2)$. Since $N=N_0+R$,
$N_0(s,t)=1/(\pi|s-t|^2)$, and $R$ is symmetric with
$\int_{\Sone}R(s,t)\,dt=r(s)$, the last two displays and
\eqref{eq:R-representation} give
\begin{equation}\label{eq:bulk-boundary-energy}
 \mathcal B(h,k)
 =\frac12\iint N(s,t)(h(s)-h(t))\cdot(k(s)-k(t))\,ds\,dt
      -\int_{\Sone}r\,h\cdot k\,ds.
\end{equation}
Finally, polarization of
Lemma~\ref{lem:curvature} and the conformal change of variables
$x=\Phi(z)$ give
\[
 E(h,k)=\mathcal B(h,k)
       +\int_{\Sone}|\Phi'|(\kappa\circ\Phi)h\cdot k\,ds.
\]
Indeed, under this change of variables the Dirichlet integral of each
component is invariant, $\Lambda_G^{\mathrm N}\,dx$ on $G$ becomes $dq$
on $\D$, and arclength on $\partial G$ becomes $|\Phi'|\,ds$ on
$\Sone$. This is \eqref{eq:boundary-energy}.

All double integrals converge absolutely:
$|h(s)-h(t)|\le C|s-t|$ and $|k(s)-k(t)|\le C|s-t|$, so the integrand of
the term with $N_0$ is bounded, and $R$ is integrable.
\end{proof}

The identities \eqref{eq:harmonic-boundary-energy} and
\eqref{eq:bulk-boundary-energy} hold for every measure $q$ as in
\eqref{eq:subcritical} and for arbitrary Lipschitz data
$h,k:\Sone\to\R^m$. Tangency and the special form
\eqref{eq:conformal-potential} of $q$ are needed only for the last step of
the proof.

\begin{lemma}[The boundary gradient of an eigenfunction; {\cite[(4.16)--(4.17) and Section~6]{OpenAI}}]
\label{lem:boundary-gradient}
Let $0\ne u\in E_1(G)$ and put
\begin{equation}
 g(s)=(\nabla u)(\Phi(s)),\qquad s\in\Sone.
 \label{eq:gradient-boundary}
\end{equation}
Then $g$ is smooth and tangent, and
\begin{equation}
 W_g=(\nabla u)\circ\Phi,
 \qquad X_g=\nabla u,
 \qquad E(g)=0.
 \label{eq:gradient-zero-energy}
\end{equation}
Moreover, $g$ is not identically zero.
\end{lemma}

\begin{proof}
The function $g$ is smooth by \eqref{eq:eigenfunction-regularity}, and
it is tangent because $\partial_\nu u=0$ on $\partial G$. Each
derivative $\partial_iu$ satisfies
$-\Delta\partial_iu=\Lambda_G^{\mathrm N}\partial_iu$ in $G$. Since
$\Delta(f\circ\Phi)=|\Phi'|^2(\Delta f)\circ\Phi$, each component of
$(\nabla u)\circ\Phi$ is smooth on $\overline\D$, satisfies
$-\Delta W=aW$ with $a=\Lambda_G^{\mathrm N}|\Phi'|^2$, and has the
corresponding component of $g$ as its trace.  Uniqueness
in Lemma~\ref{lem:resolvent}\,\textup{(iv)} gives $W_g=(\nabla u)\circ\Phi$,
hence $X_g=\nabla u$, and
Lemma~\ref{lem:energy-nullspace} gives $E(g)=0$. If $g$ were
identically zero, the same uniqueness would give $W_g=0$, hence
$\nabla u=0$ in $G$. Then $u$ would be constant, and the eigenfunction
equation would give $u=0$.
\end{proof}

\begin{proposition}[Multiplier identity; {\cite[Proposition~4.4]{OpenAI}}]\label{prop:multiplier}
Let $u$ and $g$ be as in Lemma~\ref{lem:boundary-gradient}.
For every Lipschitz function $b:\Sone\to\R$,
\begin{equation}
 E(bg)=\frac12\iint_{\Sone\times\Sone}
 N(s,t)(b(s)-b(t))^2g(s)\cdot g(t)\,ds\,dt.
 \label{eq:multiplier}
\end{equation}
\end{proposition}

\begin{proof}
For tangent Lipschitz $h$ and $\theta\in\R$, nonnegativity of $E$ and
$E(g)=0$ give
\[
 0\le E(g+\theta h)
       =2\theta E(g,h)+\theta^2E(h).
\]
Thus $E(g,h)=0$, and in particular $E(g,b^2g)=0$.
Subtract \eqref{eq:boundary-energy} for this pair from the same
formula for $(bg,bg)$. The single integrals over $\Sone$ in the two
formulas are equal, because $(bg)\cdot(bg)=g\cdot(b^2g)$.
For the double integral, the identity
\[
 |\alpha x-\beta y|^2
 -(x-y)\cdot(\alpha^2x-\beta^2y)
       =(\alpha-\beta)^2x\cdot y
\]
with $\alpha=b(s)$, $\beta=b(t)$, $x=g(s)$, and $y=g(t)$ gives
\[
 E(bg)=E(bg,bg)-E(g,b^2g)
 =\frac12\iint N(s,t)(b(s)-b(t))^2g(s)\cdot g(t)\,ds\,dt.
\]
This proves \eqref{eq:multiplier}.

Absolute convergence follows from
$|b(s)-b(t)|\le C_b|s-t|$, the boundedness of $g$, and the
integrability of $R$.
\end{proof}

\section{The reciprocal-kernel inequality}\label{sec:reciprocal}

The reciprocal-kernel argument follows \cite[Section~5]{OpenAI}.
We give an alternative proof of its Green-kernel regularization and
state the quadrature argument as a general resolvent-preservation
result. The column update and passage to the boundary follow
\cite[Lemmas~5.4 and~5.6]{OpenAI}. The normalization below explains the
connection with complete Nevanlinna--Pick kernels.

For a real symmetric matrix $A$, write $n_+(A)$ for the number of its
strictly positive eigenvalues, counted with multiplicity, and put
\[
 \mathcal P=\{A:n_+(A)\le1\}.
\]
For a kernel $k$ with values in $(0,\infty]$ and points
$x_1,\ldots,x_m$, the \emph{reciprocal matrix} of $k$ at these points
is the matrix with entries $1/k(x_i,x_j)$, where $1/\infty=0$; for a
matrix $M$ with positive entries, the reciprocal matrix of $M$ is
$(1/M_{ij})$.
The condition that reciprocal matrices belong to $\mathcal P$ appears
in the characterization of complete Nevanlinna--Pick kernels by
McCullough and Quiggin
\cite{McCullough1992,McCullough1994,Quiggin1993,Quiggin1994}; see
\cite[Sections~1 and~5]{OpenAI} and
\cite[Corollary~1.12]{AglerMcCarthy2000}.
The following operations preserve $\mathcal P$:
\begin{itemize}
\item congruence $A\mapsto V^TAV$ with an invertible matrix $V$, in
particular with a diagonal matrix with positive entries (Sylvester's
law of inertia);
\item restriction to a principal submatrix (Cauchy's interlacing
theorem);
\item subtraction of a positive semidefinite matrix (monotonicity of
the ordered eigenvalues);
\item limits of matrices of a fixed size (continuity of the
eigenvalues).
\end{itemize}

\begin{lemma}[Normalization at a distinguished point]
\label{lem:reciprocal-normalization}
Let $x_1,\ldots,x_m,z$ be points at which a symmetric kernel $k$ has
strictly positive values, allowing $+\infty$ on the diagonal.
Suppose $k(x_i,z)<\infty$ and put
\[
 B_{ij}=\frac1{k(x_i,x_j)},\qquad
 \ell_i=\frac1{k(x_i,z)},\qquad
 H_{ij}=\frac{k(x_i,z)k(z,x_j)}{k(x_i,x_j)}.
\]
If $\mu=k(z,z)<\infty$, then
\begin{equation}\label{eq:finite-pole-normalization}
 \begin{pmatrix}B&\ell\\\ell^T&\mu^{-1}\end{pmatrix}
 \in\mathcal P
 \quad\Longleftrightarrow\quad
 \left(1-\frac{H_{ij}}\mu\right)_{i,j=1}^m\succeq0.
\end{equation}
If $k(z,z)=+\infty$, then
\begin{equation}\label{eq:infinite-pole-normalization}
 \begin{pmatrix}B&\ell\\\ell^T&0\end{pmatrix}
 \in\mathcal P
 \quad\Longleftrightarrow\quad
 \alpha^TH\alpha\le0
 \quad\text{whenever }\mathbf1^T\alpha=0.
\end{equation}
\end{lemma}

\begin{proof}
In the finite case, congruence with
$\operatorname{diag}(k(x_1,z),\ldots,k(x_m,z),\mu)$ gives
\[
 \begin{pmatrix}H&\mu\mathbf1\\\mu\mathbf1^T&\mu\end{pmatrix}
 \sim
 \begin{pmatrix}H-\mu\mathbf1\mathbf1^T&0\\0&\mu\end{pmatrix}.
\]
Since $\mu>0$, this proves \eqref{eq:finite-pole-normalization}.

In the infinite case, congruence with
$\operatorname{diag}(k(x_1,z),\ldots,k(x_m,z),1)$ gives
\[
 \begin{pmatrix}H&\mathbf1\\\mathbf1^T&0\end{pmatrix}.
\]
Choose a basis of $\mathbf1^\perp$ and a vector $e$ with
$\mathbf1^Te=1$. In this basis the last matrix becomes
\[
 \begin{pmatrix}C&b&0\\b^T&a&1\\0&1&0\end{pmatrix}
 \sim C\oplus\begin{pmatrix}a&1\\1&0\end{pmatrix},
\]
where $C$ represents the restriction of $H$ to $\mathbf1^\perp$.
The second congruence follows by replacing the last coordinate $s$
by $s-b^Tu$. The last $2\times2$ block has one positive and one
negative eigenvalue. Consequently
\[
 n_+\!\begin{pmatrix}B&\ell\\\ell^T&0\end{pmatrix}
 =1+n_+(H|_{\mathbf1^\perp}),
\]
which proves \eqref{eq:infinite-pole-normalization}.
\end{proof}

The finite normalization in
\eqref{eq:finite-pole-normalization} is the one used in the
McCullough--Quiggin characterization of complete Nevanlinna--Pick
kernels. See \cite[Corollary~1.12 and Theorem~3.1]{AglerMcCarthy2000}.
The infinite-diagonal version gives the conditional negativity needed
here directly from the bordered reciprocal matrix. It is not obtained
by substituting $\mu=+\infty$ into the finite normalization.

\begin{proposition}[Reciprocal-kernel inequality; {\cite[Theorem~5.1]{OpenAI}}]\label{thm:kernel}
Let $dq=a(x)\,dx$ on $\D$, where $a\in L^\infty(\D)$,
$a\ge0$, and assume \eqref{eq:subcritical} with $d<1$.
For the kernels $K,N$ of Lemma~\ref{lem:kernels} and $p\in\D$, set
\begin{equation}\label{eq:distance-kernel}
 D_p(s,t)=\frac{K(p,s)K(p,t)}{N(s,t)}\quad(s\ne t),
 \qquad D_p(s,s)=0.
\end{equation}
Then $D_p$ is conditionally negative semidefinite. More precisely,
for every finite
family of points $s_i\in\Sone$ and coefficients $\alpha_i\in\R$,
$1\le i\le m$,
\begin{equation}\label{eq:conditional-negativity}
 \sum_{i=1}^m\alpha_i=0
 \quad\Longrightarrow\quad
 \sum_{i,j=1}^m\alpha_i\alpha_jD_p(s_i,s_j)\le0.
\end{equation}
\end{proposition}

In the rest of Section~\ref{sec:reciprocal}, $a$, $q$, and $d$ are as
in Proposition~\ref{thm:kernel}, and $T$, $S$, $K$, $R$, $N$ are the
corresponding objects of Lemmas~\ref{lem:resolvent}
and~\ref{lem:kernels}; the measure $q$ need not be the measure
\eqref{eq:conformal-potential}.

\paragraph{Plan of the proof.}
Lemma~\ref{lem:regularization} constructs bounded continuous kernels
increasing to $\mathcal G$, with reciprocal matrices in $\mathcal P$.
Lemma~\ref{lem:update} preserves this property when a positive
multiple of one kernel column times its transpose is added.
Proposition~\ref{prop:resolvent-preservation} then passes from these
finite updates to positive resolvent perturbations by quadrature.

Removing the regularization gives
Proposition~\ref{prop:interior-reciprocal}. Boundary limits give the
bordered matrix in Lemma~\ref{lem:border}, and the infinite-pole
normalization in Lemma~\ref{lem:reciprocal-normalization} gives
conditional negativity. The truncation of the potential is removed
only after the boundary limits have been taken.

\subsection{Regularization and finite updates}

\begin{lemma}[Regularized disk Green kernels; {\cite[Lemma~5.3]{OpenAI}}]\label{lem:regularization}
There are bounded, continuous, strictly positive kernels
$\mathcal G_\gamma$ on $\D\times\D$, $0<\gamma\le1$, such that
\[
 \mathcal G_\gamma\uparrow\mathcal G\quad(\gamma\downarrow0),
 \qquad
 \left(\mathcal G_\gamma(x_i,x_j)^{-1}\right)_{i,j=1}^m\in\mathcal P
\]
for every finite family $x_1,\ldots,x_m\in\D$ of points, not
necessarily distinct. On the diagonal,
$\mathcal G_\gamma(x,x)=1/(2\pi\gamma)$ and $\mathcal G(x,x)=+\infty$.
\end{lemma}

\begin{proof}
For $x,y\in\D$, put
\[
 c(x,y)=1+\frac{2|x-y|^2}{(1-|x|^2)(1-|y|^2)}.
\]
Thus $c(x,y)=\cosh d_{\D}(x,y)$ for the Poincar\'e distance
associated with the line element $2|dz|/(1-|z|^2)$. The disk Green
formula \eqref{eq:disk-green} gives, for $x\ne y$,
\[
 \mathcal G(x,y)=\frac1{2\pi}\operatorname{artanh}\frac1{c(x,y)}.
\]
Set
\[
 A_x=\frac{1+|x|^2}{1-|x|^2},\qquad
 B_x=\frac{2x}{1-|x|^2},\qquad v_x=\frac{B_x}{A_x}.
\]
Then $|v_x|<1$ and
\begin{equation}\label{eq:hyperbolic-rank-one}
 c(x,y)=A_xA_y-B_x\cdot B_y.
\end{equation}
Consequently, finite matrices of $c$ belong to $\mathcal P$.
Moreover,
\[
 \frac1{c(x,y)}
 =\frac1{A_xA_y}\sum_{j=0}^{\infty}(v_x\cdot v_y)^j
\]
is a positive semidefinite kernel. Indeed, each summand is a Gram
kernel, and the series converges on every finite family of points.
The Schur product theorem gives the same conclusion for $c^{-m}$
for every positive integer $m$.

We use the following instance of Kaluza's sign criterion
\cite[Theorem~1.3]{BariczVestiVuorinen2011}, including its short proof.
Write
\[
 H(z)=\sum_{n=0}^{\infty}a_nz^n
     =\frac{\operatorname{artanh}\sqrt z}{\sqrt z},
 \qquad a_n=\frac1{2n+1},\qquad 0\le z<1,
\]
where $H(0)=1$. Define $d_n$ by
\[
 d_n=a_n-\sum_{k=1}^{n-1}d_k a_{n-k},\qquad n\ge1.
\]
The ratios $r_n=a_n/a_{n-1}=(2n-1)/(2n+1)$ increase with $n$.
Since $d_1=a_1>0$, induction and the recursion at $n-1$ give
\[
 d_n=\sum_{k=1}^{n-1}d_k
       \bigl(r_na_{n-1-k}-a_{n-k}\bigr)\ge0,
 \qquad n\ge2.
\]
Also $d_n\le a_n$. Hence multiplication of the absolutely convergent
series yields
\[
 H(z)^{-1}=1-\sum_{n=1}^{\infty}d_nz^n.
\]
Since $H(z)\to\infty$ as $z\uparrow1$, monotone convergence gives
$\sum_{n\ge1}d_n=1$.

It follows that
\[
 F(x,y)=\sum_{n=1}^{\infty}d_n c(x,y)^{-(2n-1)}
\]
is a positive semidefinite kernel. The series converges uniformly on
$\D\times\D$, so $F$ is continuous. Furthermore,
\[
 0<F(x,y)\le1,\qquad F(x,x)=1,\qquad
 \mathcal G(x,y)^{-1}=2\pi\bigl(c(x,y)-F(x,y)\bigr).
\]
The last identity holds also on the diagonal with $1/\infty=0$.
For $0<\gamma\le1$, define
\begin{equation}\label{eq:regularized-kernel}
 \mathcal G_\gamma(x,y)^{-1}
 =2\pi\bigl(c(x,y)-(1-\gamma)F(x,y)\bigr).
\end{equation}
Its finite matrices belong to $\mathcal P$ by
\eqref{eq:hyperbolic-rank-one} and the positive semidefiniteness of $F$.
The denominator is at least $2\pi\gamma$, because $c\ge1$ and $F\le1$.
Thus $\mathcal G_\gamma$ is continuous and
\[
 0<\mathcal G_\gamma\le(2\pi\gamma)^{-1},\qquad
 \mathcal G_\gamma(x,x)=(2\pi\gamma)^{-1}.
\]
Finally $F\ge0$ implies $\mathcal G_\gamma\uparrow\mathcal G$ as
$\gamma\downarrow0$, including on the diagonal.
\end{proof}

\begin{lemma}[Column update; {\cite[Lemma~5.4]{OpenAI}}]\label{lem:update}
Let $M$ be a finite real symmetric matrix with strictly positive
entries whose reciprocal matrix $(1/M_{ij})$ belongs to $\mathcal P$.
Let $z$ be an index and $c\ge0$, and put
$\widetilde M_{ij}=M_{ij}+cM_{iz}M_{zj}$. Then
$(1/\widetilde M_{ij})\in\mathcal P$.
\end{lemma}

\begin{proof}
Write $\mu=M_{zz}>0$. Positive diagonal congruence of the reciprocal
matrix with $\operatorname{diag}(M_{iz})$ gives
\[
 H_{ij}=\frac{M_{iz}M_{zj}}{M_{ij}},
 \qquad H\in\mathcal P,
 \qquad H_{iz}=H_{zi}=\mu .
\]
Set
\[
 U:=\mu\mathbf1\mathbf1^T-H.
\]
Lemma~\ref{lem:reciprocal-normalization} shows that the principal
submatrix with index $z$ removed is positive semidefinite.
The $z$ row and column vanish, so $U\succeq0$.
Since $H_{ii}>0$,
\[
 0\le U_{ii}<\mu,
 \qquad |U_{ij}|\le\sqrt{U_{ii}U_{jj}}<\mu.
\]
The same diagonal congruence transforms the updated reciprocal into
$H_{ij}/(1+cH_{ij})$. Entrywise,
\[
 \begin{aligned}
 \frac \mu{1+c\mu}-\frac{H_{ij}}{1+cH_{ij}}
 &=\frac{U_{ij}}{(1+c\mu)(1+c\mu-cU_{ij})}\\
 &=\sum_{n\ge0}
       \frac{c^nU_{ij}^{n+1}}{(1+c\mu)^{n+2}}.
 \end{aligned}
\]
The series converges absolutely. Its coefficients are nonnegative,
and every Schur power of $U$ is positive semidefinite. Thus the
updated reciprocal is congruent to a positive rank-one matrix minus
a positive semidefinite matrix, and belongs to $\mathcal P$.
\end{proof}

\subsection{Positive resolvent perturbations}

\begin{proposition}[Preservation under a positive resolvent perturbation]
\label{prop:resolvent-preservation}
Let $X$ be a compact metric space and let
$k\colon X\times X\to(0,\infty)$ be continuous and symmetric.
Suppose every finite reciprocal matrix of $k$ belongs to $\mathcal P$.
Let $\rho$ be a finite nonnegative Borel measure on $X$, and suppose
the integral operator
\[
 (T_\rho f)(x)=\int_X k(x,y)f(y)\,d\rho(y)
 \quad\text{on }L^2(\rho)
\]
satisfies $\|T_\rho\|<1$. Then
\begin{equation}\label{eq:general-positive-resolvent}
 k^\rho(x,y)=k(x,y)+
 \big\langle k(x,\cdot),(I-T_\rho)^{-1}k(y,\cdot)
                                      \big\rangle_{L^2(\rho)}
\end{equation}
is continuous, symmetric, and strictly positive. Every finite
reciprocal matrix of $k^\rho$ belongs to $\mathcal P$.
\end{proposition}

The proof is the positive-quadrature argument in
\cite[Proposition~5.5]{OpenAI}, with its kernel-independent content
stated separately. No additional positive-semidefiniteness assumption
on $k$ is required. In fact, this property follows from the reciprocal
condition when all kernel values are finite. To see this, fix $z\in X$
and write
\[
 F(x,y)=1-\frac{k(x,z)k(z,y)}{k(x,y)k(z,z)},
 \qquad \delta(x)=\frac{k(x,z)}{\sqrt{k(z,z)}}.
\]
Lemma~\ref{lem:reciprocal-normalization} gives $F\succeq0$ and
$0\le F(x,x)<1$. Hence $|F(x,y)|<1$, and the Schur product theorem gives
\[
 k(x,y)=\delta(x)\delta(y)\sum_{n=0}^\infty F(x,y)^n\succeq0
\]
on each finite family. The argument below uses only the reciprocal
condition and nonnegativity of the kernel chains.

\begin{proof}
If $\rho=0$, the assertion is immediate. Otherwise choose finite Borel
partitions of $X$ with cell diameters at most $\varepsilon$.
Discard cells of zero mass, choose a node $z_\alpha$ in each remaining
cell, and put $w_\alpha=\rho(\text{cell}_\alpha)>0$.
Let $\pi_\varepsilon$ map each such cell to its node, and let
$T_\varepsilon$ have kernel
$k(\pi_\varepsilon x,\pi_\varepsilon y)$ on $L^2(\rho)$.
Uniform continuity gives
\[
 \begin{aligned}
 \|T_\varepsilon-T_\rho\|
 &\le\rho(X)\,
   \operatorname*{ess\,sup}_{x,y}
   |k(\pi_\varepsilon x,\pi_\varepsilon y)-k(x,y)|\longrightarrow0,\\
 \sup_{x\in X}\|k(x,\pi_\varepsilon(\cdot))-k(x,\cdot)\|_{L^2(\rho)}
 &\longrightarrow0.
 \end{aligned}
\]
Thus $\|T_\varepsilon\|<1$ for small $\varepsilon$, and the resolvent
identity gives convergence in operator norm of
$(I-T_\varepsilon)^{-1}$ to $(I-T_\rho)^{-1}$.
Consequently
\begin{equation}\label{eq:quadrature-limit}
 \begin{aligned}
 k_\varepsilon^\rho(x,y)
 &:={k}(x,y)+
 \big\langle k(x,\pi_\varepsilon(\cdot)),
 (I-T_\varepsilon)^{-1}k(y,\pi_\varepsilon(\cdot))\big\rangle\\
 &\longrightarrow k^\rho(x,y)
 \end{aligned}
\end{equation}
uniformly on $X\times X$.

In the orthonormal basis of normalized cell indicators, set
\begin{equation}\label{eq:quadrature-matrix}
 \mathsf T_{\alpha\beta}
  =\sqrt{w_\alpha}\,k(z_\alpha,z_\beta)\sqrt{w_\beta},
 \qquad
 (\xi_x)_\alpha=\sqrt{w_\alpha}\,k(x,z_\alpha).
\end{equation}
Then $\|\mathsf T\|=\|T_\varepsilon\|<1$ and
\[
 k_\varepsilon^\rho(x,y)
 =k(x,y)+\xi_x^T(I-\mathsf T)^{-1}\xi_y.
\]
Fix a finite family of evaluation points, and append the set $Z$ of
nodes to obtain an indexed family $\mathcal J$. Repetitions cause no
problem. Duplicating points replaces a matrix $A$ by a pullback
$V^TAV$, which cannot increase the number $n_+$ of positive eigenvalues.
For $Y\subset Z$ and $i,j\in\mathcal J$, define
\begin{equation}\label{eq:chain-sum}
 M^Y_{ij}=k(i,j)+
 \sum_{r\ge1}\ \sum_{y_1,\ldots,y_r\in Y}
 k(i,y_1)w_{y_1}k(y_1,y_2)\cdots w_{y_r}k(y_r,j).
\end{equation}
The nonnegative series converges, since
\[
 0<M^Y_{ij}\le M^Z_{ij}
 =k(i,j)+\xi_i^T(I-\mathsf T)^{-1}\xi_j<\infty.
\]
For $z\in Z\setminus Y$, group the chains according to their number
of visits to $z$. Each segment between consecutive visits contains
only nodes of $Y$. Therefore
\[
 M^{Y\cup\{z\}}_{ij}
 =M^Y_{ij}+
 \sum_{r\ge1}w_z^rM^Y_{iz}(M^Y_{zz})^{r-1}M^Y_{zj}.
\]
Finiteness and strict positivity imply $w_zM^Y_{zz}<1$, so
\begin{equation}\label{eq:resolvent-node-update}
 M^{Y\cup\{z\}}_{ij}
 =M^Y_{ij}+\frac{w_z}{1-w_zM^Y_{zz}}M^Y_{iz}M^Y_{zj}.
\end{equation}
Starting with $M^\varnothing_{ij}=k(i,j)$, apply
Lemma~\ref{lem:update} one node at a time. The reciprocal matrix of
$M^Z$ belongs to $\mathcal P$. Its principal submatrix on the evaluation
points is the reciprocal matrix of $k_\varepsilon^\rho$.
Now use \eqref{eq:quadrature-limit} and the closedness of $\mathcal P$.

Symmetry and continuity of $k^\rho$ follow from
\eqref{eq:general-positive-resolvent}. The Neumann series consists of
nonnegative kernel chains, so $k^\rho\ge k>0$.
\end{proof}

Let $0\le\rho\le q$ have compact support in $\D$. Denote by
$T_\rho$ the Green operator on $L^2(\rho)$ and set, for $x\ne y$,
\begin{equation}\label{eq:interior-resolvent}
 M^\rho(x,y)=\mathcal G(x,y)
 +\left\langle\mathcal G(x,\cdot),
       (I-T_\rho)^{-1}\mathcal G(y,\cdot)\right\rangle_{L^2(\rho)}.
\end{equation}
The second term is finite because $\mathcal G(x,\cdot)\in L^2(\rho)$ and
$\|(I-T_\rho)^{-1}\|\le(1-d)^{-1}$ by
Lemma~\ref{lem:resolvent}\,\textup{(ii)} and~\textup{(v)}.
Set $M^\rho(x,x)=+\infty$. The kernel $M^\rho$ is the Green kernel of
the Dirichlet operator $-\Delta-a_\rho$ on $\D$, where $a_\rho$ is the
density of $\rho$. This identification
is not used in the proofs below; it is explained, with its
probabilistic description, in Section~\ref{subsec:prob-occupation}.

\begin{proposition}[Reciprocal property of the Green resolvent; {\cite[Proposition~5.5]{OpenAI}}]\label{prop:interior-reciprocal}
Let $\rho$ be a measure with $0\le\rho\le q$ and compact support in
$\D$, and let $x_1,\ldots,x_m\in\D$ be distinct. Then the matrix
\[
 A_{ij}=\begin{cases}
 M^\rho(x_i,x_j)^{-1},&i\ne j,\\
 0,&i=j,
 \end{cases}
\]
belongs to $\mathcal P$.
\end{proposition}

\begin{proof}
Fix $\gamma\in(0,1]$ and let $C=\mathcal G_\gamma$.
Choose a compact metric subspace $L\Subset\D$ containing
$\operatorname{supp}\rho$ and all the evaluation points.
The integral operator $T_C$ on $L^2(\rho)$ satisfies
\[
 |T_Cf|\le T_\rho|f|,\qquad \|T_C\|\le d<1,
\]
because $0<C\le\mathcal G$.
Lemma~\ref{lem:regularization} and
Proposition~\ref{prop:resolvent-preservation}, applied on $L$, show
that the reciprocal matrices of
\begin{equation}\label{eq:regularized-resolvent}
 M^{\rho,\gamma}(x,y)=\mathcal G_\gamma(x,y)+
 \big\langle\mathcal G_\gamma(x,\cdot),
 (I-T_C)^{-1}\mathcal G_\gamma(y,\cdot)\big\rangle_{L^2(\rho)}
\end{equation}
belong to $\mathcal P$.

Expanding the resolvent gives
\[
 \begin{aligned}
 M^{\rho,\gamma}(x,y)
 &=\mathcal G_\gamma(x,y)\\
 &\quad+\sum_{r\ge1}\int_{\D^r}
 \mathcal G_\gamma(x,y_1)\mathcal G_\gamma(y_1,y_2)\cdots
 \mathcal G_\gamma(y_r,y)
 \,d\rho(y_1)\cdots d\rho(y_r).
 \end{aligned}
\]
Every integrand is nonnegative and increases as $\gamma\downarrow0$.
Monotone convergence in each integral and in the series yields
\[
 M^{\rho,\gamma}(x,y)\uparrow M^\rho(x,y)
 \qquad(x\ne y).
\]
The limit is finite because its resolvent term is at most
\[
 \frac{\|\mathcal G(x,\cdot)\|_{L^2(\rho)}
       \|\mathcal G(y,\cdot)\|_{L^2(\rho)}}{1-d}.
\]
On the diagonal,
$M^{\rho,\gamma}(x,x)\ge\mathcal G_\gamma(x,x)\to+\infty$.
Thus the reciprocal matrices converge to the matrix in the
proposition, with zero diagonal. Closedness of $\mathcal P$ completes
the proof.
\end{proof}

\subsection{Boundary limits and conditional negativity}

\begin{lemma}[Bordered boundary matrix; {\cite[Lemma~5.6]{OpenAI}}]\label{lem:border}
For distinct $s_1,\ldots,s_m\in\Sone$ and $p\in\D$, put
\[
 B_{ij}=\begin{cases}N(s_i,s_j)^{-1},&i\ne j,\\0,&i=j,\end{cases}
 \qquad \ell_i=K(p,s_i)^{-1}.
\]
Then
\begin{equation}\label{eq:bordered-matrix}
 \begin{pmatrix}B&\ell\\\ell^T&0\end{pmatrix}\in\mathcal P.
\end{equation}
\end{lemma}

\begin{proof}
First fix compactly supported $0\le\rho\le q$. A superscript
$\rho$ denotes the kernels obtained with this measure. Apply
Proposition~\ref{prop:interior-reciprocal} to
$rs_1,\ldots,rs_m,p$, where $r>|p|$ and $r\uparrow1$.
The explicit disk Green function gives
\begin{equation}\label{eq:green-boundary-limits}
 \frac{\mathcal G(rs,x)}{1-r}\longrightarrow P_s(x),
 \qquad
 \frac{\mathcal G(rs,rt)}{(1-r)^2}
       \longrightarrow\frac1{\pi|s-t|^2}=N_0(s,t)\quad(s\ne t).
\end{equation}
The first convergence is uniform for $x$ in compact subsets of $\D$.
The second follows from
\[
 \mathcal G(rs,rt)
 =\frac1{4\pi}\log\left(1+
       \frac{(1-r^2)^2}{r^2|s-t|^2}\right).
\]
Since the support of $\rho$ is compact and
$\|(I-T_\rho)^{-1}\|\le(1-d)^{-1}$,
\[
 \begin{aligned}
 \frac{M^\rho(rs_i,rs_j)}{(1-r)^2}
 &\longrightarrow N_0(s_i,s_j)+
   \langle P_{s_i},(I-T_\rho)^{-1}P_{s_j}\rangle_{L^2(\rho)}
   =N^\rho(s_i,s_j),\\
 \frac{M^\rho(p,rs_i)}{1-r}
 &\longrightarrow P_{s_i}(p)+
   \langle\mathcal G(p,\cdot),(I-T_\rho)^{-1}P_{s_i}
                  \rangle_{L^2(\rho)}
   =K^\rho(p,s_i).
 \end{aligned}
\]
The functions $P_{s_i}$ are bounded on $\operatorname{supp}\rho$, so
the inner products are finite. Congruence with
$\operatorname{diag}(1-r,\ldots,1-r,1)$ transforms the reciprocal
matrix of Proposition~\ref{prop:interior-reciprocal} into the matrix
with entries $(1-r)^2/M^\rho(rs_i,rs_j)$ for $i\ne j$, with last row
and column $(1-r)/M^\rho(p,rs_i)$, and with zero diagonal. Its limit as
$r\uparrow1$ is \eqref{eq:bordered-matrix} for
$K^\rho,N^\rho$. Hence that
bordered matrix belongs to $\mathcal P$.

Now let $d\rho_n=\mathbf1_{\{|x|\le r_n\}}\,dq$, with $r_n\uparrow1$.
The series representations \eqref{eq:kernel-chains}, whose terms are
nonnegative, and monotone convergence give
\[
 K^{\rho_n}(p,s)\uparrow K(p,s),\qquad
 N^{\rho_n}(s,t)\uparrow N(s,t)\quad(s\ne t).
\]
Their limits are finite and positive by Lemma~\ref{lem:kernels}.
Taking the entrywise reciprocals and passing to the limit proves
\eqref{eq:bordered-matrix}. The argument does not use
$P_s\in L^2(q)$, which may fail.
\end{proof}

\begin{proof}[Proof of Proposition~\ref{thm:kernel}]
For distinct $s_1,\ldots,s_m\in\Sone$, let $B$ and $\ell$ be as in
Lemma~\ref{lem:border}. Congruence of its bordered matrix with
\[
 \operatorname{diag}(K(p,s_1),\ldots,K(p,s_m),1)
\]
gives
\[
 \begin{pmatrix}
 (D_p(s_i,s_j))_{i,j=1}^m&\mathbf1\\
 \mathbf1^T&0
 \end{pmatrix}\in\mathcal P.
\]
The infinite-pole argument in
Lemma~\ref{lem:reciprocal-normalization} proves
\eqref{eq:conditional-negativity}. Repeated points are handled by
adding their coefficients.
\end{proof}

\section{Proof of Theorem~\ref{thm:smooth}}
\label{sec:hilbert-finish}

Throughout this section, $q$ is the measure
\eqref{eq:conformal-potential}, and $K$, $N$, and $D_p$ are the
corresponding kernels. The last step of the proof uses a theorem
recalled by Schoenberg \cite[Section~3, pp.~525--526]{Schoenberg}: a
conditionally negative semidefinite kernel that vanishes on the
diagonal is the squared distance of a map into a Hilbert space. We
follow \cite[Section~6]{OpenAI} for this construction, the summed
multiplier identity, and the final case distinction. We verify that
the map is Lipschitz, as required by Proposition~\ref{prop:multiplier}.
We use the classical smooth-domain multiplicity bound recalled in
Lemma~\ref{lem:smooth-multiplicity}. No result of the companion
Lipschitz paper is needed.

\begin{lemma}[A Lipschitz map into Hilbert space with squared distance $D_p$; {\cite[Lemma~6.1]{OpenAI}}]
\label{lem:hilbert-realization}
For each $p\in\D$ there are a separable real Hilbert space $\mathcal H$
and an injective Lipschitz map $b:\Sone\to\mathcal H$ such that
\begin{equation}\label{eq:hilbert-distance}
 \|b(s)-b(t)\|_{\mathcal H}^2
 =D_p(s,t)=\frac{K(p,s)K(p,t)}{N(s,t)},\qquad s\ne t.
\end{equation}
After an isometric embedding $\mathcal H\hookrightarrow\ell^2$, write
$b=(b_1,b_2,\ldots)$. Every $b_j$ is a Lipschitz function, and
\begin{equation}\label{eq:coordinate-distance}
 \sum_{j=1}^{\infty}|b_j(s)-b_j(t)|^2=D_p(s,t).
\end{equation}
\end{lemma}

\begin{proof}
Fix $s_0\in\Sone$ and set
\[
 k_p(s,t)=\tfrac12\bigl[D_p(s,s_0)+D_p(t,s_0)-D_p(s,t)\bigr],
 \qquad s,t\in\Sone.
\]
This kernel is positive semidefinite. Indeed, for points
$s_1,\ldots,s_m\in\Sone$ and real numbers $c_1,\ldots,c_m$, put
$c_0=-\sum_{i=1}^mc_i$. Since $D_p(s_0,s_0)=0$,
Proposition~\ref{thm:kernel}, applied to $s_0,\ldots,s_m$ and
$c_0,\ldots,c_m$, gives
\[
 \sum_{i,j=1}^mc_ic_jk_p(s_i,s_j)
 =-\frac12\sum_{i,j=0}^mc_ic_jD_p(s_i,s_j)\ge0.
\]
Let $V$ be the real vector space of finitely supported functions on
$\Sone$, and let $\delta_s\in V$ be the indicator function of $\{s\}$.
The kernel $k_p$ defines a positive semidefinite bilinear form on $V$
with $\langle\delta_s,\delta_t\rangle=k_p(s,t)$. Let $\mathcal H_0$ be
the completion of the quotient of $V$ by the null space of this form,
and let $b(s)\in\mathcal H_0$ be the class of $\delta_s$. Then
$\langle b(s),b(t)\rangle=k_p(s,t)$, and
\[
 \|b(s)-b(t)\|^2=k_p(s,s)+k_p(t,t)-2k_p(s,t)=D_p(s,t).
\]
Lemma~\ref{lem:kernels} gives
$\|K(p,\cdot)\|_\infty<\infty$, and $N\ge N_0$. Hence
\begin{equation}\label{eq:b-lipschitz}
 0<D_p(s,t)\le
 \pi\|K(p,\cdot)\|_\infty^2|s-t|^2
 \qquad(s\ne t).
\end{equation}
Thus $b$ is injective and Lipschitz. The closed span $\mathcal H$ of
its image is separable, because $b$ is continuous. We embed
$\mathcal H$ isometrically into $\ell^2$.
Coordinate projections and Parseval's identity give the remaining
assertions.
\end{proof}

\begin{proposition}[The summed multiplier identity; {\cite[(6.2)]{OpenAI}}]
\label{prop:summed-multiplier}
Let $0\ne u\in E_1(G)$, let $g=(\nabla u)\circ\Phi$ on $\Sone$ be as in
Lemma~\ref{lem:boundary-gradient},
and let $b=(b_j)_{j\ge1}$ be the realization in
Lemma~\ref{lem:hilbert-realization} associated with $p\in\D$.
Then
\begin{equation}\label{eq:sum-energy}
 \sum_{j=1}^{\infty}E(b_jg)
 =\frac12\left|\int_{\Sone}K(p,s)g(s)\,ds\right|^2
 =\frac12|\nabla u(\Phi(p))|^2.
\end{equation}
Every summand is nonnegative.
\end{proposition}

\begin{proof}
By Proposition~\ref{prop:multiplier},
\[
 \sum_{j=1}^mE(b_jg)
 =\frac12\iint N(s,t)\sum_{j=1}^m|b_j(s)-b_j(t)|^2
                  g(s)\cdot g(t)\,ds\,dt.
\]
Every summand on the left is nonnegative by
Lemma~\ref{lem:energy-nullspace}, because $b_jg$ is tangent and
Lipschitz.
Off the diagonal,
\[
 0\le N(s,t)\sum_{j=1}^m|b_j(s)-b_j(t)|^2
 \le N(s,t)D_p(s,t)=K(p,s)K(p,t).
\]
Since $K(p,\cdot)$ and $g$ are bounded, dominated convergence and
\eqref{eq:coordinate-distance} give
\[
 \sum_{j\ge1}E(b_jg)
 =\frac12\iint K(p,s)K(p,t)g(s)\cdot g(t)\,ds\,dt
 =\frac12|W_g(p)|^2,
\]
where the last equality is \eqref{eq:K-representation}.
Equation~\eqref{eq:gradient-zero-energy} completes the proof.
\end{proof}

\begin{proof}[Proof of Theorem~\ref{thm:smooth}]
Fix $0\ne u\in E_1(G)$ and put $g=(\nabla u)\circ\Phi$ on $\Sone$.
By Lemma~\ref{lem:boundary-gradient}, the set
\begin{equation}\label{eq:nonempty-Sg}
 S_g=\{s\in\Sone:g(s)\ne0\}
\end{equation}
is nonempty and open. Lemma~\ref{lem:smooth-multiplicity} gives
$\dim E_1(G)\le2$.

Suppose that $\nabla u(y_0)=0$ for some $y_0\in G$. Set
$p=\Phi^{-1}(y_0)$ and choose the injective Lipschitz map
$b=(b_j)_{j\ge1}$ from Lemma~\ref{lem:hilbert-realization}.
Proposition~\ref{prop:summed-multiplier} and
Lemma~\ref{lem:energy-nullspace} imply
\[
 0=\sum_{j\ge1}E(b_jg),\qquad
 E(b_jg)=0,\qquad X_{b_jg}=\nabla v_j,
 \qquad v_j\in E_1(G).
\]
Here every summand is nonnegative. Elliptic regularity gives
$v_j\in C^\infty(\overline G)$, and
Lemma~\ref{lem:resolvent}\,\textup{(iv)} makes $X_{b_jg}$ continuous on
$\overline G$ with the prescribed boundary values. Thus
\begin{equation}\label{eq:boundary-multiplier-gradients}
 \nabla v_j(\Phi(s))=b_j(s)g(s)
 \qquad(s\in\Sone,\ j\ge1).
\end{equation}
We complete the argument by the two cases in
\cite[Section~6]{OpenAI}.

Suppose first that $S_g$ is dense in $\Sone$. The set $b(S_g)$ is not
contained in an affine line. Otherwise continuity would put
$b(\Sone)$ in that line. A continuous injection of the circle into a
line is impossible, since its compact connected image would be a
nondegenerate interval, and removing an interior point disconnects an
interval but not a circle with one point removed.

We may therefore choose $s_1,s_2,s_3\in S_g$ such that
$b(s_2)-b(s_1)$ and $b(s_3)-b(s_1)$ are linearly independent. Some two
coordinates $j,k$ then satisfy
\[
 \det\begin{pmatrix}
  1&b_j(s_1)&b_k(s_1)\\
  1&b_j(s_2)&b_k(s_2)\\
  1&b_j(s_3)&b_k(s_3)
 \end{pmatrix}\ne0.
\]
Indeed, if every two-coordinate minor vanished, any nonzero
coordinate of the first difference vector would force the second
difference vector to be proportional to it. If
$\alpha u+\beta v_j+\gamma v_k=0$, taking boundary gradients at these
three points gives
\[
 \bigl(\alpha+\beta b_j(s_i)+\gamma b_k(s_i)\bigr)g(s_i)=0
 \qquad(i=1,2,3).
\]
Since $g(s_i)\ne0$, the determinant implies
$\alpha=\beta=\gamma=0$. Thus $u,v_j,v_k$ are linearly independent,
contrary to $\dim E_1(G)\le2$.

Suppose now that $S_g$ is not dense. There is a nonempty open arc
$I\subset\Sone$ on which $g=0$. Injectivity of $b$ and openness of
$S_g$ give an index $j$ for which $b_j$ is nonconstant on $S_g$.
Equation~\eqref{eq:boundary-multiplier-gradients} shows that $u$ and
$v_j$ are linearly independent. Their gradients vanish on the
connected boundary arc $\Gamma=\Phi(I)$, so
\[
 u|_\Gamma=c_0,\qquad v_j|_\Gamma=c_1
\]
for constants $c_0,c_1$. Choose $(\alpha,\beta)\ne(0,0)$ such that
$\alpha c_0+\beta c_1=0$. The nonzero eigenfunction
$w=\alpha u+\beta v_j$ then satisfies
\[
 (\Delta+\Lambda_G^{\mathrm N})w=0\quad\hbox{in }G,
 \qquad w=\partial_\nu w=0\quad\hbox{on }\Gamma.
\]

Choose a ball $B$ centered at an interior point of $\Gamma$ so small
that $B\cap\partial G\subset\Gamma$, and extend $w$ by zero to
$B\setminus\overline G$. Write $\widetilde w$ for this extension.
For every $\varphi\in C_c^\infty(B)$, Green's identity gives
\[
 \int_B\widetilde w
       (\Delta+\Lambda_G^{\mathrm N})\varphi\,dx
 =\int_{B\cap\partial G}
       \bigl(w\partial_\nu\varphi-\varphi\partial_\nu w\bigr)\,d\sigma
 =0.
\]
This is a local version of the zero-extension argument used in
\cite[Lemma, p.~414]{Filonov}. Interior elliptic regularity implies
that $\widetilde w$ is smooth and real analytic in $B$. It vanishes
on the nonempty exterior open set $B\setminus\overline G$, hence on
$B$. Analytic continuation in the connected domain $G$ now gives
$w=0$ in $G$, a contradiction. Only smoothness of the boundary was
used in this argument.

Both cases are impossible. Therefore $\nabla u\ne0$ in $G$.
The function $u$ is continuous on the compact set $\overline G$.
An interior maximum or minimum would be a critical point, so
\[
 \min_{\partial G}u<u(x)<\max_{\partial G}u
 \qquad(x\in G).
\]
This proves (HS1).
\end{proof}

\paragraph{A quantitative consequence for a simple eigenvalue.}
\label{subsec:quantitative-gradient}

The summed multiplier identity also gives a lower bound for the
gradient. The following formulation separates the spectral gap from
the boundary quantities that would have to be controlled in a limiting
argument.

\begin{proposition}[A quantitative gradient bound]
\label{prop:quantitative-gradient}
Suppose that $G$ is a bounded simply connected smooth planar domain
and that $\Lambda=\Lambda_G^{\mathrm N}$ is simple. Write
$\Lambda_2=\Lambda_{G,2}^{\mathrm N}$ for the next positive Neumann
eigenvalue, counted with multiplicity. Let $0\ne u\in E_1(G)$,
$g=(\nabla u)\circ\Phi$ on $\Sone$, and
\begin{equation}\label{eq:boundary-gradient-measure}
 d\mu(s)=|g(s)|^2|\Phi'(s)|\,ds,
 \qquad M=\mu(\Sone)>0.
\end{equation}
There is a constant $C_G<\infty$ such that
\begin{equation}\label{eq:trace-energy-constant}
 \|\operatorname{Tr}Y\|_{L^2(\partial G)}^2
 \le C_G\Q(Y),\qquad Y\in\mathcal X(G).
\end{equation}
For any such constant and every $p\in\D$,
\begin{equation}\label{eq:quantitative-gradient}
 |\nabla u(\Phi(p))|^2
 \ge \frac{1-\Lambda/\Lambda_2}{C_GM}
 \iint_{\Sone\times\Sone}D_p(s,t)\,d\mu(s)\,d\mu(t)>0.
\end{equation}
If $G$ is convex and
$B(x_*,r)\subset G\subset B(x_*,R)$, one may take
\begin{equation}\label{eq:convex-trace-constant}
 C_G=\frac2r\bigl(\Lambda^{-1}+R\Lambda^{-1/2}\bigr).
\end{equation}
\end{proposition}

\begin{proof}
The trace theorem, the curvature identity
\eqref{eq:curvature-energy}, and
\eqref{eq:vector-poincare} imply
\eqref{eq:trace-energy-constant}. For completeness, put
$K_- =\|\max\{-\kappa,0\}\|_{L^\infty(\partial G)}$.
The trace inequality with a small parameter gives
\[
 \|\operatorname{Tr}Y\|_2^2
 \le\varepsilon\|\nabla Y\|_2^2+C_\varepsilon\|Y\|_2^2,
 \qquad
 \Q(Y)\ge\|\nabla Y\|_2^2-K_-\|\operatorname{Tr}Y\|_2^2.
\]
Choose $\varepsilon$ with $K_-\varepsilon\le1/2$ and use
$\Lambda\|Y\|_2^2\le\Q(Y)$. It follows that
$\|Y\|_{H^1(G)}^2\le C\Q(Y)$, and the trace theorem proves
\eqref{eq:trace-energy-constant}.

The positive spectrum of the divergence--curl form is the union of
the positive Neumann and Dirichlet spectra, with multiplicities
\cite[Theorem~3.4]{Rohleder2025}. Moreover,
\[
 \Lambda_2<\lambda_G^{\mathrm D}
\]
by \cite[Corollary~4.3]{Rohleder2025}. Consequently, simplicity gives
\begin{equation}\label{eq:orthogonal-vector-gap}
 Y\in\mathcal X(G),\quad
 \langle Y,\nabla u\rangle_{L^2(G)}=0
 \quad\Longrightarrow\quad
 \Q(Y)\ge\Lambda_2\|Y\|_{L^2(G)}^2.
\end{equation}

Fix $p$ and choose $b=(b_j)_{j\ge1}$ as in
Lemma~\ref{lem:hilbert-realization}. Set
\[
 c_j=\frac{\langle X_{b_jg},\nabla u\rangle_{L^2(G)}}
              {\|\nabla u\|_{L^2(G)}^2},
 \qquad
 Y_j=X_{b_jg}-c_j\nabla u=X_{(b_j-c_j)g}.
\]
Then $Y_j\perp\nabla u$. Since $E(g,k)=0$ for every tangent
Lipschitz $k$, equations \eqref{eq:energy-definition} and
\eqref{eq:orthogonal-vector-gap} give
\[
 E(b_jg)=E((b_j-c_j)g)
 =\Q(Y_j)-\Lambda\|Y_j\|_2^2
 \ge(1-\Lambda/\Lambda_2)\Q(Y_j).
\]
Applying \eqref{eq:trace-energy-constant} to $Y_j$ yields
\begin{equation}\label{eq:coordinate-variance-bound}
 E(b_jg)\ge
 \frac{1-\Lambda/\Lambda_2}{C_G}
 \int_{\Sone}|b_j-c_j|^2\,d\mu.
\end{equation}
For every real $c$, the weighted variance identity gives
\[
 \int|b_j-c|^2\,d\mu
 \ge\frac1{2M}\iint|b_j(s)-b_j(t)|^2\,d\mu(s)\,d\mu(t).
\]
Sum \eqref{eq:coordinate-variance-bound} for $1\le j\le m$ and let
$m\to\infty$. Monotone convergence,
\eqref{eq:coordinate-distance}, and \eqref{eq:sum-energy} give
\[
 \frac12|\nabla u(\Phi(p))|^2
 \ge\frac{1-\Lambda/\Lambda_2}{2C_GM}
       \iint D_p(s,t)\,d\mu(s)\,d\mu(t).
\]
This is \eqref{eq:quantitative-gradient}. The measure $\mu$ is
nonatomic and has positive mass, because $g$ is nonzero on a nonempty
open arc. Since $D_p(s,t)>0$ for $s\ne t$, the double integral is
strictly positive.

Finally, suppose that $G$ is convex. Supporting lines and the
inclusion $B(x_*,r)\subset G$ give
$(x-x_*)\cdot\nu(x)\ge r$ on $\partial G$. The divergence theorem
and Cauchy--Schwarz yield
\[
 \begin{aligned}
 r\|\operatorname{Tr}Y\|_{L^2(\partial G)}^2
 &\le\int_{\partial G}(x-x_*)\cdot\nu\,|Y|^2\,ds\\
 &=2\|Y\|_{L^2(G)}^2
    +2\int_G\sum_{i,j}(x_i-x_{*,i})Y_j\partial_iY_j\,dx\\
 &\le2\|Y\|_{L^2(G)}^2
      +2R\|Y\|_{L^2(G)}\|\nabla Y\|_{L^2(G)}.
 \end{aligned}
\]
The identity follows first for smooth fields and then for $H^1$
fields by density and continuity of the trace. Convexity gives
$\kappa\ge0$, so
\[
 \|\nabla Y\|_2^2\le\Q(Y),
 \qquad \|Y\|_2^2\le\Lambda^{-1}\Q(Y).
\]
Substitution proves \eqref{eq:convex-trace-constant}.
\end{proof}

\begin{remark}[Uniformity for convex domains]
\label{rem:convex-quantitative-uniformity}
For smooth convex domains satisfying
$B(x_*,r)\subset G\subset B(x_*,R)$ with fixed $r,R$, the
trace constant can be chosen uniformly. Indeed, the
Payne--Weinberger inequality \cite{PayneWeinberger1960} gives
$\Lambda\ge\pi^2/\operatorname{diam}(G)^2\ge\pi^2/(4R^2)$.
Thus \eqref{eq:convex-trace-constant} gives
\begin{equation}\label{eq:convex-uniform-trace-constant}
 C_G\le C_*:=\frac{4R^2}{r}
                 \left(\frac2{\pi^2}+\frac1\pi\right).
\end{equation}
No upper bound on curvature is required. In particular, this estimate
does not deteriorate solely because a smooth convex approximation
develops corners. For $\|u\|_{L^2(G)}=1$, it also gives
\[
 M=\|\operatorname{Tr}\nabla u\|_{L^2(\partial G)}^2
 \le C_*\Q(\nabla u)=C_*\Lambda^2.
\]

This proves uniform control of the trace factor, but does not give a
uniform positive lower bound for the full right-hand side of
\eqref{eq:quantitative-gradient}. Such a conclusion also requires
control of $1-\Lambda/\Lambda_2$ and of
\[
 \frac1M\iint D_p(s,t)\,d\mu(s)\,d\mu(t).
\]
The spectral factor tends to zero if simple eigenvalues approach a
double eigenvalue. Even with a uniform spectral gap, the present
kernel bounds do not supply a uniform lower bound for this weighted
boundary variance. An upper bound for $M$ alone does not exclude
concentration of $\mu$. We therefore make no uniform gradient claim
for arbitrary convex approximating families.
\end{remark}

\begin{remark}[Where the hypotheses are used]
\label{rem:hypotheses}
The proof of Theorem~\ref{thm:smooth} uses the multiplicity bound
$\dim E_1(G)\le2$, rather than simplicity. Apart from this multiplicity
bound, the first-eigenvalue hypothesis enters through two analytic
facts. The first is the inequality
\eqref{eq:strict-spectral-gap}, which gives \eqref{eq:subcritical} with
$d=\Lambda_G^{\mathrm N}/\lambda_G^{\mathrm D}<1$ for the measure
\eqref{eq:conformal-potential}. The second is
Proposition~\ref{prop:vector}, which gives the nonnegativity of $E$ and
the description of its null space in Lemma~\ref{lem:energy-nullspace};
these are used in Lemma~\ref{lem:boundary-gradient},
Proposition~\ref{prop:multiplier},
Proposition~\ref{prop:summed-multiplier}, and the proof of
Theorem~\ref{thm:smooth}.

Let $\Lambda'$ be a Neumann eigenvalue of $G$ with
$\Lambda_G^{\mathrm N}<\Lambda'<\lambda_G^{\mathrm D}$, let $u$ be a
corresponding eigenfunction, and define $q$, the kernels, $W_h$, $X_h$,
and the form $E$ with $\Lambda'$ in place of $\Lambda_G^{\mathrm N}$.
Then \eqref{eq:subcritical} holds with
$d=\Lambda'/\lambda_G^{\mathrm D}<1$, so
Lemmas~\ref{lem:resolvent}, \ref{lem:kernels},
and~\ref{lem:boundary-energy} and Proposition~\ref{thm:kernel} apply.
Put $g=(\nabla u)\circ\Phi$ on $\Sone$. The uniqueness statement of
Lemma~\ref{lem:resolvent}\,\textup{(iv)} gives $X_g=\nabla u$ as before.
For every tangent Lipschitz map $h$ we have $E(g,h)=0$. Indeed,
$\operatorname{curl}\nabla u=0$, $\Delta u=-\Lambda'u$, and
$\operatorname{Tr}X_h\cdot\nu=0$, so
\[
 E(g,h)=\int_G\Delta u\,\operatorname{div}X_h\,dx
        -\Lambda'\int_G\nabla u\cdot X_h\,dx
 =-\Lambda'\int_G\operatorname{div}(uX_h)\,dx=0.
\]
Consequently the multiplier identity \eqref{eq:multiplier} holds, with
the same proof. The identity \eqref{eq:sum-energy} also holds, because
its proof by dominated convergence does not use the sign of the terms
$E(b_jg)$. Hence $\sum_jE(b_jg)=0$ if $\nabla u(\Phi(p))=0$. The proof
cannot be continued from this point: the form $E$ is no longer
nonnegative, so the vanishing of the sum does not imply the vanishing
of each term. To see that $E$ takes negative values, let
$0\ne u_1\in E_1(G)$ and $g_1=(\nabla u_1)\circ\Phi$. Write
$F_{\Lambda'}(Y)=\Q(Y)-\Lambda'\|Y\|_2^2$.
For a field $Y\in\mathcal X(G)$ with the same trace as $X_h$, put
$Z=Y-X_h\in H^1_0(G;\R^2)$. Lemma~\ref{lem:curvature} and
$-\Delta X_h=\Lambda'X_h$ give
\[
 F_{\Lambda'}(Y)-F_{\Lambda'}(X_h)
 =\int_G|\nabla Z|^2\,dx-\Lambda'\int_G|Z|^2\,dx
 \ge(\lambda_G^{\mathrm D}-\Lambda')\|Z\|_2^2\ge0.
\]
Taking $h=g_1$ and $Y=\nabla u_1$ gives
\[
 E(g_1)\le\Q(\nabla u_1)-\Lambda'\int_G|\nabla u_1|^2\,dx
 =(\Lambda_G^{\mathrm N}-\Lambda')\int_G|\nabla u_1|^2\,dx<0.
\]
\end{remark}

\begin{remark}[Topology and dimension]
\label{rem:topology-dimension}
Simple connectedness is needed already in
Proposition~\ref{prop:vector}. To see the obstruction, let $G$ be a
bounded smooth planar domain with at least two boundary components.
Let $\eta$ be harmonic in $G$, equal to one on one boundary component and
zero on the others. Then $\eta$ is nonconstant and
\[
 X=\nabla^\perp\eta\in\mathcal X(G),\qquad
 \operatorname{div}X=\operatorname{curl}X=0,\qquad X\ne0.
\]
Thus $\Q(X)=0<\Lambda_G^{\mathrm N}\|X\|_2^2$. The harmonic field
that was shown to vanish in the decomposition
\eqref{eq:vector-helmholtz} need not vanish on a domain with holes.

In dimension $n\ge3$, the reciprocal-kernel inequality already fails
for zero potential on the unit ball $\mathbb B_n$. Write
$\sigma_{n-1}=|\mathbb S^{n-1}|$. The Poisson kernel and its inward
boundary derivative give
\[
 P_t(x)=\frac{1-|x|^2}{\sigma_{n-1}|x-t|^n},\qquad
 K(0,s)=\frac1{\sigma_{n-1}},\qquad
 N(s,t)=\frac2{\sigma_{n-1}|s-t|^n}.
\]
Consequently the proposed squared-distance kernel is
\begin{equation}\label{eq:ball-dimensional-obstruction}
 D_0^{(n)}(s,t)=\frac{|s-t|^n}{2\sigma_{n-1}}.
\end{equation}
If it were conditionally negative semidefinite, the construction in
Lemma~\ref{lem:hilbert-realization} would give a Hilbert-space map $b$
with
\[
 \|b(s)-b(t)\|=(2\sigma_{n-1})^{-1/2}|s-t|^{n/2}.
\]
Join $s$ to $t$ by a spherical arc of length $L$ and divide the arc
into $m$ equal pieces. The triangle inequality gives
\[
 \|b(s)-b(t)\|
 \le (2\sigma_{n-1})^{-1/2}L^{n/2}m^{1-n/2}
 \longrightarrow0.
\]
This contradicts \eqref{eq:ball-dimensional-obstruction} for $s\ne t$.
The exponent is exactly one in the plane, where the squared chordal
distance does have the required representation.
\end{remark}

\section{Acknowledgments}

The author gratefully acknowledges Krzysztof Burdzy for many valuable discussions on hot spots over the years, particularly in recent months.  The author thanks Greg Buzzard for many conversations on the topic of this paper and for obtaining a review of the manuscript from Claude (Anthropic), which offered interesting insights and raised further questions. He also thanks Greg and Chris Janjigian for many conversations on the use of AI in the mathematical sciences and STEM fields more broadly.

\section{AI disclosure}
This manuscript was prepared with assistance from ChatGPT (OpenAI)
in exposition, calculations, and proof checking.
Theorem~\ref{thm:smooth}, the principal arguments, and the proof
structure are due to OpenAI \cite{OpenAI}. This paper provides a detailed
exposition of that proof. The sources of the individual statements
and arguments are identified in the text. The alternative regularization
proof, the general resolvent formulation, the quantitative consequences,
and the probabilistic identifications are described in the Introduction,
along with their classical inputs and their relation to the original
proof. The author does not claim that this exposition constitutes
an independent proof of OpenAI's theorem.

At present, the author does not intend to submit this manuscript
for journal publication.

%\section*{Acknowledgments}
%The author thanks Greg Buzzard for many conversations on the topic of
%this paper and for obtaining a review of the manuscript from Claude
%(Anthropic), which offered interesting insights and raised further
%questions.
%
%\section*{AI disclosure}
%This manuscript was prepared with assistance from ChatGPT (OpenAI)
%in exposition, calculations, and proof checking.
%Theorem~\ref{thm:smooth}, the principal arguments, and the proof
%structure are due to \cite{OpenAI}. This paper gives a detailed
%exposition following that proof. The sources of the individual
%statements and arguments are identified in the text. The alternative regularization
%proof, the general resolvent formulation, the quantitative consequences,
%and the probabilistic identifications are the contributions described
%in the Introduction, with their classical inputs and their relation
%to the original proof. 
%%The author
%%takes responsibility for the mathematical content of this paper.

\appendix
\section{Probabilistic interpretation}
\label{sec:probability}

The probabilistic interpretations presented in this appendix are
supplementary.  None  of its results are needed for the proof of the
main theorem, nor do they provide a complete alternative proof.

The coupling and survival-probability arguments in
\cite{BB,Pascu2002,AB,BP2004,BPP} motivate seeking a stochastic route
to the smooth-domain theorem. We begin by identifying the kernels
in its proof with weighted exit distributions, occupation densities,
and boundary jump intensities. The gauge transformation relates these
objects to a common diffusion, and the weighted Douglas identity
expresses the boundary energy through its trace process.

These identifications use classical stochastic analysis. Their purpose
here is to formulate the comparison that a probabilistic proof would
still need. Remark~\ref{rem:naim} identifies that question as conditional
negative semidefiniteness of the reciprocal Na\"im kernel. The
Feynman--Kac and weighted Douglas identities hold for bounded
nonnegative potentials satisfying \eqref{eq:subcritical}. The pointwise
Na\"im boundary limit is stated for smooth potentials. We conclude
with the nodal-domain viewpoint and the approximation constructions
of \cite{LipschitzCompanion}.

All Brownian motions have generator $\Delta$.
The symbols $\mathbb P_x$ and $\mathbb E_x$ refer to a process started
at $x$.

\subsection{A Feynman--Kac formula for \texorpdfstring{$W_h$}{W(h)}}
\label{subsec:prob-feynman-kac}

Let $a$, $q$, and $d$ be as in \eqref{eq:subcritical}, and let $T$,
$S$, $Ah$, and $W_h$ be as in Lemma~\ref{lem:resolvent}. The function
$w=W_1=1+TS1$ satisfies $w\ge1$, because $A1=1$ and $T$ and $S$
preserve positivity.
Let $(B_t)$ be planar Brownian motion, and put
\[
 \tau_\D=\inf\{t\ge0:B_t\notin\D\},\qquad
 \mathcal A_t=\int_0^{t\wedge\tau_\D}a(B_s)\,ds.
\]

\begin{proposition}[Feynman--Kac formula for $W_h$]
\label{prop:feynman-kac}
For $f\in L^2(q)$ and Lipschitz boundary data $h:\Sone\to\R^m$,
\begin{equation}\label{eq:green-occupation}
 Tf(x)=\mathbb E_x\!\left[\int_0^{\tau_\D} a(B_s)f(B_s)\,ds\right],
 \qquad
 Ah(x)=\mathbb E_x[h(B_{\tau_\D})],
\end{equation}
and
\begin{equation}\label{eq:stopped-feynman-kac}
 W_h(x)=\mathbb E_x\!\left[e^{\mathcal A_{\tau_\D}}h(B_{\tau_\D})\right],
 \qquad x\in\D.
\end{equation}
The expectation in \eqref{eq:stopped-feynman-kac} is absolutely
convergent, and
\begin{equation}\label{eq:finite-gauge}
 \sup_{x\in\D}\mathbb E_xe^{\mathcal A_{\tau_\D}}
 \le\|w\|_\infty
 \le1+\frac{C_aq(\D)^{1/2}}{1-d}<\infty.
\end{equation}
The same holds when $q$ is replaced by a measure $\rho$ with
$0\le\rho\le q$, with $a$ replaced by the density $a_\rho$ of $\rho$ in
$\mathcal A_t$, $T$, $S$, $W_h$, $w$, and $C_a$, and with the same
constant $d$.
\end{proposition}

These identities are classical. See
\cite[Chapters~2 and~4]{ChungZhao1995}. We give a short proof that uses
\eqref{eq:subcritical} and $w\ge1$.
With Brownian generator $\frac12\Delta$, a factor $\frac12$ must be
inserted in the occupation integral and in $\mathcal A_t$.

\begin{proof}
The identities \eqref{eq:green-occupation} are the classical
representations of the Green operator and of the Poisson integral of
the disk. By \eqref{eq:green-evaluation-bound}, the occupation integral
converges absolutely at every $x$.

The identity $W_h=Ah+TW_h$ and the strong Markov property imply
that
\[
 M_t=W_h(B_{t\wedge\tau_\D})+
       \int_0^{t\wedge\tau_\D}a(B_s)W_h(B_s)\,ds
\]
is a continuous martingale. It is bounded on each finite time
interval. Integration by parts with $e^{\mathcal A_t}$ cancels
the finite-variation term. Since $a$ and $W_h$ are bounded,
the resulting stochastic integral is a true martingale on each
finite interval. Consequently,
\begin{equation}\label{eq:finite-time-gauge}
 W_h(x)=\mathbb E_x\!\left[
 e^{\mathcal A_t}W_h(B_{t\wedge\tau_\D})\right].
\end{equation}

Taking $h=1$ and using $w\ge1$ gives
$\mathbb E_x e^{\mathcal A_t}\le w(x)$. Monotone convergence yields
\eqref{eq:finite-gauge}; the bound on $\|w\|_\infty$ follows from
$w=1+TS1$, \eqref{eq:green-evaluation-bound}, and
\eqref{eq:L2-resolvent-bound}.
Dominated convergence in \eqref{eq:finite-time-gauge}, with dominator
$\|W_h\|_\infty e^{\mathcal A_{\tau_\D}}$, proves
\eqref{eq:stopped-feynman-kac}. The same proof applies componentwise
to vector data. If $0\le\rho\le q$, its density is bounded by $a$
and its form bound has the same $d$.
\end{proof}

\subsection{\texorpdfstring{$K$}{K} as a weighted exit density}
\label{subsec:prob-exit-law}

We keep the notation of Section~\ref{subsec:prob-feynman-kac}. The
kernel $K$ of \eqref{eq:kernel-K} is the density, with respect to
arclength, of the exit distribution of $B$ from $\D$ weighted by
$e^{\mathcal A_{\tau_\D}}$.
Indeed, Proposition~\ref{prop:feynman-kac} and
Lemma~\ref{lem:kernels}, applied to
boundary test functions, give the equality of finite Borel measures
\begin{equation}\label{eq:weighted-exit-law}
 K(p,s)\,ds
 =\mathbb E_p\!\left[e^{\mathcal A_{\tau_\D}};\ B_{\tau_\D}\in ds\right],
 \qquad
 \int_{\Sone}K(p,s)\,ds=w(p).
\end{equation}
Thus $K(p,s)/w(p)$ is a probability density. More precisely, the
probability law on stopped Brownian paths defined by
\[
 \frac{d\mathbb Q_p}{d\mathbb P_p}
       =\frac{e^{\mathcal A_{\tau_\D}}}{w(p)}
\]
satisfies
\[
 \frac{W_h(p)}{w(p)}=\mathbb E_p^{\mathbb Q}[h(B_{\tau_\D})].
\]
The Markov property gives its density up to time $t\wedge\tau_\D$ as
$e^{\mathcal A_t}w(B_{t\wedge\tau_\D})/w(p)$. Since
$\Delta w=-aw$, its interior generator is
\begin{equation}\label{eq:gauge-diffusion}
 L^wf=w^{-1}(\Delta+a)(wf)
      =\Delta f+2\nabla\log w\cdot\nabla f.
\end{equation}
The quotient $W_h/w$ is used again in
Section~\ref{subsec:prob-douglas}.

\subsection{The measure \texorpdfstring{$q$}{q} and the exit time from \texorpdfstring{$G$}{G}}
\label{subsec:prob-clock}

Let $G$ and $\Phi$ be as in Section~\ref{sec:preliminaries}, let $q$ be
the measure \eqref{eq:conformal-potential}, and let $\mathcal A_t$ and
$\tau_\D$ be as in Section~\ref{subsec:prob-feynman-kac} with
$a=\Lambda_G^{\mathrm N}|\Phi'|^2$. Then $\mathcal A_t$ is
$\Lambda_G^{\mathrm N}$ times the time change that turns $B$ into
Brownian motion in $G$. Set
\[
 C_t=\int_0^{t\wedge\tau_\D}|\Phi'(B_s)|^2\,ds,
 \qquad \eta_t=\inf\{s\ge0:C_s>t\}.
\]
Up to its lifetime, $\widehat B_t=\Phi(B_{\eta_t})$ is Brownian
motion in $G$. If $\tau_G$ denotes its exit
time, set $\widehat B_{\tau_G}=\Phi(B_{\tau_\D})$ by continuity.
Conformal invariance gives
\begin{equation}\label{eq:conformal-clock}
 \tau_G=C_{\tau_\D},\qquad
 \mathcal A_{\tau_\D}=\Lambda_G^{\mathrm N}\tau_G.
\end{equation}
Consequently \eqref{eq:stopped-feynman-kac} becomes
\begin{equation}\label{eq:physical-exit-law}
 W_h(p)=\mathbb E_{\Phi(p)}\!\left[
       e^{\Lambda_G^{\mathrm N}\tau_G}h\bigl(\Phi^{-1}(\widehat B_{\tau_G})\bigr)
                                      \right].
\end{equation}
For a measure $0\le\rho\le q$ with density $a_\rho\le a$, as used in the
truncations of Section~\ref{sec:reciprocal}, the corresponding additive
functional satisfies
$\int_0^{\tau_\D}a_\rho(B_s)\,ds\le\mathcal A_{\tau_\D}
=\Lambda_G^{\mathrm N}\tau_G$.

For $0\ne u\in E_1(G)$ and $g=(\nabla u)\circ\Phi$ as in
Lemma~\ref{lem:boundary-gradient}, formula
\eqref{eq:physical-exit-law} reads, in the domain $G$,
\begin{equation}\label{eq:gradient-exit-law}
 \nabla u(y)=\mathbb E_y\!\left[
                  e^{\Lambda_G^{\mathrm N}\tau_G}\nabla u(\widehat B_{\tau_G})\right],
 \qquad y\in G.
\end{equation}
This is a weighted average of the boundary gradient, with total mass
$w(\Phi^{-1}(y))$.

\subsection{A weighted Douglas identity}
\label{subsec:prob-douglas}

Let $a$, $q$, and $d$ be as in \eqref{eq:subcritical}, and let $W_h$,
$K$, $R$, and $N$ be as in Lemmas~\ref{lem:resolvent}
and~\ref{lem:kernels}. Let $w=W_1$ be as in
Section~\ref{subsec:prob-feynman-kac}, and put
$r(s)=\int_{\Sone}R(s,t)\,dt$ and $H_h=W_h/w$.

\begin{proposition}[Weighted Douglas identity]
\label{prop:weighted-douglas}
For scalar or vector Lipschitz boundary data $h,k$,
\begin{equation}\label{eq:weighted-Douglas}
 \int_\D w^2\nabla H_h:\nabla H_k\,dx
 =\frac12\iint_{\Sone\times\Sone}N(s,t)
       (h(s)-h(t))\cdot(k(s)-k(t))\,ds\,dt.
\end{equation}
\end{proposition}

\begin{proof}
The function $w-1=Tw$ has zero trace, and boundedness of $a$ and $w$
gives $w\in W^{2,p}(\D)$ for every finite $p$, so its normal derivative
on $\Sone$ is well defined and continuous. Let $h$ be a Lipschitz
function on $\Sone$. Since $w-1\in H^1_0(\D)$ and $Ah$ is harmonic and
belongs to $H^1(\D)$, we have $\int_\D\nabla w\cdot\nabla Ah\,dx=0$.
The divergence theorem for $Ah\,\nabla w\in W^{1,1}(\D)$ and
$\Delta w=-aw$ therefore give
\[
 \int_{\Sone}h\,\partial_\nu w\,ds=\int_\D Ah\,\Delta w\,dx
 =-\int_{\Sone}h(s)\int_\D P_s(x)w(x)\,dq(x)\,ds.
\]
Since $h$ is arbitrary, $-\partial_\nu w(s)=\int_\D P_sw\,dq$ for almost
every $s$. Inserting $w(x)=\int_{\Sone}K(x,t)\,dt$, which is
\eqref{eq:K-representation} with $h=1$, we obtain from Tonelli's
theorem and the definition \eqref{eq:kernel-R} of $R$
\begin{equation}\label{eq:gauge-boundary-flux}
 -\partial_\nu w(s)
 =\int_\D P_s(x)w(x)\,dq(x)
 =\int_{\Sone}R(s,t)\,dt=r(s)
 \qquad\text{for almost every }s\in\Sone.
\end{equation}
In particular $r\in L^\infty(\Sone)$.

Let $\mathcal B(h,k)$ be the form defined in the proof of
Lemma~\ref{lem:boundary-energy}; formula
\eqref{eq:bulk-boundary-energy} holds for scalar or vector Lipschitz
data under \eqref{eq:subcritical}. Now $w\ge1$, $w=1$ on $\Sone$, and
$\Delta w=-aw$. Expanding $W_h=wH_h$
and $W_k=wH_k$ and integrating the cross terms by means of
$\operatorname{div}(w\nabla w)=|\nabla w|^2-aw^2$ give
\[
 \mathcal B(h,k)
 =\int_\D w^2\nabla H_h:\nabla H_k\,dx
                  +\int_{\Sone}\partial_\nu w\,h\cdot k\,ds.
\]
Equations~\eqref{eq:gauge-boundary-flux} and
\eqref{eq:bulk-boundary-energy} prove
\eqref{eq:weighted-Douglas}.
\end{proof}

Formula~\eqref{eq:weighted-Douglas} has an interpretation through a
boundary trace process. On the state space $\overline\D$, with reference measure
$w^2dx$, consider the regular reflecting form
\[
 \mathcal E^w(f,f)=\int_\D w^2|\nabla f|^2\,dx,
 \qquad D(\mathcal E^w)=H^1(\D).
\]
Since $\operatorname{div}(w^2\nabla H_h)=w(\Delta+a)W_h=0$, the function
$H_h$ is the $\mathcal E^w$-harmonic extension of $h$.
The bounds $1\le w\le C$ give finite reference mass, and the energy
measures $2w^2|\nabla f|^2dx$ charge no boundary mass. Arclength is a
smooth measure with full quasi-support on $\Sone$.
Time-change its associated diffusion by the boundary additive
functional whose Revuz measure is $ds$. The resulting boundary
trace has the energy in \eqref{eq:weighted-Douglas}
\cite[Theorem~6.2]{FHY2004}.
Its form domain is $H^{1/2}(\Sone)$, and the identity extends from
Lipschitz data to this domain by density. Indeed, the weighted bulk
norms are equivalent to the unweighted ones, and
$r\in L^\infty(\Sone)$ bounds the $R$ contribution on $L^2(\Sone)$.

In the L\'evy-kernel convention, its jump kernel is $N(s,t)\,dt$.
The corresponding symmetric jump measure is
$\tfrac12N(s,t)\,ds\,dt$ when the jump-energy integral is written
without a prefactor. The identity \eqref{eq:weighted-Douglas} was
proved directly above and requires no excursion estimates.

\begin{remark}[The Na\"im normalization]
\label{rem:naim}
Suppose here that $a\in C^\infty(\overline\D)$ satisfies
\eqref{eq:subcritical}. Let $M$ be the Dirichlet Green kernel of
$-\Delta-a$, given by \eqref{eq:interior-resolvent} with $\rho=q$.
For $s\ne t$, boundary differentiation of the Green kernel gives
\begin{equation}\label{eq:naim-normalization}
 \Theta_p(s,t)
 :=\lim_{r\uparrow1}
 \frac{M(rs,rt)}{M(p,rs)M(p,rt)}
 =\frac{N(s,t)}{K(p,s)K(p,t)}
 =\frac1{D_p(s,t)}.
\end{equation}
Indeed, the numerator divided by $(1-r)^2$ tends to $N(s,t)$,
and each denominator factor divided by $1-r$ tends to its Poisson
kernel. The same limits for compactly supported bounded potentials
were proved in Lemma~\ref{lem:border}.
Thus $D_p$ is the reciprocal of the boundary Green kernel with
Na\"im's normalization \cite{Naim1957}. For this normalization and
its boundary extension in a general diffusion setting, see
\cite[(1.6) and Definition~3.15]{KajinoMurugan2026}.

The Markov interpretation has a fixed scalar normalization.
For the gauge diffusion \eqref{eq:gauge-diffusion} killed on $\Sone$,
the Green kernel relative to $w^2dx$ is
\[
 G^w(x,y)=\frac{M(x,y)}{w(x)w(y)}.
\]
Its Na\"im kernel with pole $p$ is therefore
$\Theta_p^w=w(p)^2\Theta_p$. If
$d\mu_p(s)=K(p,s)\,ds/w(p)$ is its exit law, then
\[
 \Theta_p^w(s,t)\,d\mu_p(s)\,d\mu_p(t)
 =N(s,t)\,ds\,dt.
\]
This identifies the normalization in the weighted Douglas formula.

When $a=0$, let $\phi_p(z)=(z-p)/(1-\overline pz)$. Directly from
\eqref{eq:disk-poisson} and \eqref{eq:kernel-N},
\[
 D_p(s,t)
 =\frac{(1-|p|^2)^2|s-t|^2}
        {4\pi|s-p|^2|t-p|^2}
 =\frac1{4\pi}|\phi_p(s)-\phi_p(t)|^2.
\]
The conditional negative semidefiniteness for nonzero $a$ is the
additional conclusion of Proposition~\ref{thm:kernel}. The ball
calculation in Remark~\ref{rem:topology-dimension} shows why the
existence of a symmetric Markov process alone does not imply it.
\end{remark}

\subsection{\texorpdfstring{$M^\rho$}{M-rho} as an occupation density}
\label{subsec:prob-occupation}

Let $a$, $q$, and $d$ be as in \eqref{eq:subcritical}, let $\rho$ be a
measure with $0\le\rho\le q$ and compact support in $\D$, let $a_\rho$
be its density, and let $M^\rho$ be the kernel
\eqref{eq:interior-resolvent}. If $f$ is bounded and nonnegative, then
\begin{equation}\label{eq:resolvent-occupation}
 \int_\D M^\rho(x,y)f(y)\,dy
 =\mathbb E_x\!\left[\int_0^{\tau_\D}
       \exp\left(\int_0^t a_\rho(B_s)\,ds\right)f(B_t)\,dt\right].
\end{equation}
Indeed, the form bound gives
\[
 \int_\D|\nabla\psi|^2-\int_\D a_\rho\psi^2
 \ge(1-d)\int_\D|\nabla\psi|^2,
 \qquad \psi\in H^1_0(\D).
\]
Thus the Dirichlet operator $-\Delta-a_\rho$ is strictly positive.
Integrating its killed Feynman--Kac semigroup gives its inverse
\cite[Chapter~4]{ChungZhao1995}, whose resolvent identity is
\eqref{eq:interior-resolvent}. This proves
\eqref{eq:resolvent-occupation}, using the continuous interior versions
on both sides. Elliptic regularity gives continuity of the inverse,
and Feynman--Kac smoothing together with the positive spectral gap
gives continuity of its potential representation.
Proposition~\ref{prop:interior-reciprocal} does not follow from the
positivity of $M^\rho$. Its proof uses the explicit formula for
$\mathcal G$.

\subsection{Reflecting Brownian motion killed on the nodal arc}
\label{rem:nodal-survival}

The eigenvalue problem on a nodal domain, with the Dirichlet condition
on the nodal line and the Neumann condition on $\partial\Omega$, is
studied in
Ba\~nuelos and Pang \cite{BP2004} and by Ba\~nuelos, Pang, and Pascu
\cite{BPP}.

Let $\Omega$ be a bounded simply connected Lipschitz domain. For
$0\ne u\in E_1(\Omega)$, put
$\Gamma=\{u=0\}\cap\overline\Omega$,
$\Omega_+=\{u>0\}\cap\Omega$, and
$O_+=\{u>0\}\cap\overline\Omega$.
If $X$ is reflecting Brownian motion in $\overline\Omega$ and
$\tau_\Gamma=\inf\{t\ge0:X_t\in\Gamma\}$, then
\begin{equation}\label{eq:nodal-killed-semigroup}
 P_t^+f(x)=\mathbb E_x[f(X_t);\ t<\tau_\Gamma]
\end{equation}
is the semigroup of $X$ killed at $\tau_\Gamma$ and restricted to
$L^2(\Omega_+)$. The process reflects on $\partial\Omega$ before
it is killed. In the terminology of
\cite[Section~4.4]{FOT2011}, $(P_t^+)$ is the semigroup of the part of
the Neumann form on the relatively open set $O_+\subset\overline\Omega$.
More precisely, its form domain, with functions extended by zero
outside $O_+$, is
\[
 \mathcal F_+=\{f\in H^1(\Omega):
   \widetilde f=0\ \text{quasi-everywhere on }
                     \overline\Omega\setminus O_+\}.
\]
Here $\widetilde f$ is the quasi-continuous representative for the
Neumann form, and quasi-everywhere means outside a set of zero
capacity for that form. Since $\mathcal F_+\subset H^1(\Omega)$ and
$H^1(\Omega)\hookrightarrow L^2(\Omega)$ is compact, the operators
$P_t^+$ are compact. For $f\ge0$, $P_t^+f$ is bounded below by the
corresponding expectation for Brownian motion killed on
$\partial\Omega_+$; since $\Omega_+$ is connected, the semigroup
$(P_t^+)$ is therefore positivity improving.

The function $u_+=u\mathbf1_{\Omega_+}$ belongs to $H^1(\Omega)$ by
the nodal-restriction argument in Lemma~\ref{lem:smooth-multiplicity},
which applies equally to a continuous Lipschitz-domain eigenfunction,
and its continuous representative vanishes on
$\overline\Omega\setminus O_+$. Thus $u_+\in\mathcal F_+$.
The weak Neumann equation gives
\[
 \int_\Omega\nabla u_+\cdot\nabla v\,dx
 =\Lambda_\Omega^{\mathrm N}\int_\Omega u_+v\,dx,
 \qquad v\in\mathcal F_+,
\]
because $v$ and $\nabla v$ vanish almost everywhere off $\Omega_+$.
Consequently $P_t^+u_+=e^{-\Lambda_\Omega^{\mathrm N}t}u_+$.
Since $u_+>0$ on
$\Omega_+$ and the semigroup $(P_t^+)$ is positivity improving,
$\Lambda_\Omega^{\mathrm N}$ is the lowest eigenvalue of the associated
nonnegative operator $A_+$, where $P_t^+=e^{-tA_+}$, and this eigenvalue
is simple. Equivalently, the generator has largest eigenvalue
$-\Lambda_\Omega^{\mathrm N}$. The spectral theorem and interior
parabolic smoothing give
\begin{equation}\label{eq:nodal-survival-asymptotic}
 e^{\Lambda_\Omega^{\mathrm N}t}\mathbb P_x(\tau_\Gamma>t)
 \longrightarrow
 u(x)\frac{\int_{\Omega_+}u}{\int_{\Omega_+}u^2},
\end{equation}
in $L^2(\Omega_+)$ and locally uniformly in $\Omega_+$.
This statement concerns the semigroup $(P_t^+)$ on the fixed set $O_+$;
it does not say that $\Lambda_\Omega^{\mathrm N}$ is a simple Neumann
eigenvalue of $\Omega$.

The nodal-arc theorem for Lipschitz domains in
\cite[Lemma~2.3]{LipschitzCompanion} shows that $\Omega_+$ is a Jordan domain
bounded by $\Gamma$ and an arc of $\partial\Omega$.
Let $\D_+$ be the upper half-disk and choose a conformal map
$F_+:\D_+\to\Omega_+$ carrying its diameter
onto $\Gamma$. Let $B$ be Brownian motion in $\D_+$ that is reflected
on the semicircle and killed on the diameter, at time $\sigma$. The time
change by $|F_+'|^2$ gives
\begin{equation}\label{eq:nodal-conformal-clock}
 \mathbb P_{F_+(z)}(\tau_\Gamma>t)
 =\mathbb P_z\left\{\int_0^\sigma|F_+'(B_s)|^2\,ds>t\right\},
 \qquad z\in\D_+.
\end{equation}
To identify the process, put $I=[-1,1]$ and $\rho=|F_+'|^2$.
The form of $B$ is
\[
 e(v,v)=\int_{\D_+}|\nabla v|^2\,dz,
 \qquad
 \mathcal V=\{v\in H^1(\D_+):\operatorname{Tr}_I v=0\}.
\]
Its time change has the same energy, reference measure $\rho\,dz$,
and domain $\mathcal V\cap L^2(\rho\,dz)$
\cite[Section~6.2]{FOT2011}. This is the conformal pullback of the part
form on $\mathcal F_+$. Indeed, for $v=f\circ F_+$,
\[
 \int_\Omega |f|^2\,dx=\int_{\D_+}|v|^2\rho\,dz,
 \qquad
 \int_\Omega|\nabla f|^2\,dx=\int_{\D_+}|\nabla v|^2\,dz.
\]
The domain identification follows first for bounded continuous
form-core functions supported away from $\Gamma$. The Jordan
extension of $F_+$ carries $I$ onto $\Gamma$.
Conversely, if $v$ is smooth on the closed half-disk and vanishes
near $I$, then $v\circ F_+^{-1}$ has a zero extension across $\Gamma$
belonging to $\mathcal F_+$.
Approximation on the half-disk and closure in the displayed
energy-plus-$L^2$ norm give the full domains. The two endpoints are
polar and impose no additional condition. Thus no Lipschitz
regularity of the nodal domain is being assumed.

The time-change theorem now identifies the processes outside an
exceptional set of zero capacity, and hence gives
\eqref{eq:nodal-conformal-clock} there. The identity extends to every
$z\in\D_+$ as follows.
We use the versions of $X$ and $B$ that are defined for every starting
point; started at an interior point, each is an ordinary Brownian
motion up to its first exit from a disk with closure in the interior.
Take a disk $V$ centered at $z$ with closure in $\D_+$, and apply the
strong Markov property to $B$ at the exit time of $V$ and to $X$ at
the exit time of $F_+(V)$. Up to these exit times the time change is
the classical conformal invariance of Brownian motion, which holds for
every starting point. The exit distribution of $B$ from $V$ is
normalized arclength on $\partial V$, and a subset of $\partial V$ of
zero capacity has zero arclength. Hence the exit point lies almost
surely outside the exceptional set, and
\eqref{eq:nodal-conformal-clock} at $z$ follows from
\eqref{eq:nodal-conformal-clock} at the exit point.

The paper \cite{BP2004} assumes that
$r\mapsto r^2|F_+'(re^{i\theta})|^2$ is nondecreasing.
The nodal-arc theorem alone does not give this monotonicity.

By \eqref{eq:nodal-survival-asymptotic}, the exponential moment of
$\tau_\Gamma$ at the rate $\Lambda_\Omega^{\mathrm N}$ is infinite. Let $\tau_{\partial\Omega}$ denote the first exit time
of ordinary Brownian motion from $\Omega$. Then
\begin{equation}\label{eq:nodal-critical-gauge}
 \mathbb E_xe^{\Lambda_\Omega^{\mathrm N}\tau_\Gamma}=\infty
 \quad(x\in\Omega_+),\qquad
 \sup_{x\in\Omega}\mathbb E_xe^{\alpha\tau_{\partial\Omega}}<\infty
 \quad(0\le\alpha<\lambda_\Omega^{\mathrm D}),
\end{equation}
where the second assertion is well known for bounded Lipschitz
domains. See, for example, \cite{Banuelos1991}.
Since $\Lambda_\Omega^{\mathrm N}<\lambda_\Omega^{\mathrm D}$ by
\eqref{eq:strict-spectral-gap}, the second assertion applies with
$\alpha=\Lambda_\Omega^{\mathrm N}$.

\subsection{Connections with the Lipschitz approximation}
\label{subsec:prob-approximation}

The companion paper \cite{LipschitzCompanion} uses a thin-strip
perturbation to select a prescribed first eigendirection, followed by
a conformal exhaustion. We record the probabilistic meaning of the
two constructions. These observations have no role in the proof of
Theorem~\ref{thm:smooth}.

\paragraph{The strip coefficient.}
\label{rem:strip-return-time}
Let $\Omega$ be a bounded simply connected Lipschitz domain, and
choose $h>0$ so that $\Lambda_\Omega^{\mathrm N}<\pi^2/(4h^2)$.
The coefficient in the strip-splitting formula of
\cite[Lemma~3.9]{LipschitzCompanion} is
$d_h=\sqrt{\Lambda_\Omega^{\mathrm N}}
       \tan(\sqrt{\Lambda_\Omega^{\mathrm N}}h)$.
Let $Y$ have generator $d^2/dx^2$ on $(0,h]$, with reflection at
$h$, and let $\sigma_0$ be its first hitting time of $0$.
We stop $Y$ there, with $Y_{\sigma_0}=0$. For
$0<\Lambda_\Omega^{\mathrm N}<\pi^2/(4h^2)$, the solution $g$ of
$-g''=\Lambda_\Omega^{\mathrm N}g$, $g(0)=1$, $g'(h)=0$, satisfies
\begin{equation}\label{eq:strip-probabilistic-profile}
 g(x)=\mathbb E_x e^{\Lambda_\Omega^{\mathrm N}\sigma_0}
     =\frac{\cos(\sqrt{\Lambda_\Omega^{\mathrm N}}(h-x))}{\cos(\sqrt{\Lambda_\Omega^{\mathrm N}} h)}.
\end{equation}
To verify the expectation, apply It\^o's formula to the stopped
process
$e^{\Lambda_\Omega^{\mathrm N}(t\wedge\sigma_0)}g(Y_{t\wedge\sigma_0})$.
The local-time term at $h$ vanishes because $g'(h)=0$.
Since $g\ge1$ is bounded, monotone convergence first gives
$\mathbb E_x e^{\Lambda_\Omega^{\mathrm N}\sigma_0}\le g(x)$.
Dominated convergence then gives equality. Consequently,
\begin{equation}\label{eq:strip-return-coefficient}
 d_h=g'(0)
 =\lim_{x\downarrow0}
       \frac{\mathbb E_x(e^{\Lambda_\Omega^{\mathrm N}\sigma_0}-1)}{x}
 =\sqrt{\Lambda_\Omega^{\mathrm N}}\tan(\sqrt{\Lambda_\Omega^{\mathrm N}} h).
\end{equation}
If $(u_1,u_2)$ is an orthonormal basis of a double first eigenspace,
an attachment point $p$ is called admissible when
\[
 \int_{\Omega\cap B(p,r)}|\nabla v|^2=o(r)
 \qquad(r\downarrow0,\ v\in E_1(\Omega)).
\]
For a strip of width $\varepsilon$ attached at such a point, the
matrix of the Neumann form on the two lowest positive eigenspaces,
in a suitable orthonormal basis, satisfies
\cite[Lemma~3.9]{LipschitzCompanion}
\[
 H_\varepsilon=\Lambda_\Omega^{\mathrm N}I
 -\varepsilon d_h\,\mathbf u(p)\mathbf u(p)^T+o(\varepsilon),
 \qquad \mathbf u(p)=(u_1(p),u_2(p))^T.
\]
The coefficient $-d_h$ is the number
$\int_0^h(g'^2-\Lambda_\Omega^{\mathrm N}g^2)\,dx=[gg']_0^h$. By
\eqref{eq:strip-return-coefficient}, it is expressed through the
hitting time $\sigma_0$.

\paragraph{Conformal exhaustion as a time change.}
\label{subsec:prob-timechange}
Let $F:\D\to\Omega$ be conformal, and set
\[
 F_r(z)=F(rz),\qquad G_r=F(r\D),\qquad
 \rho_r(z)=|F_r'(z)|^2=r^2|F'(rz)|^2,
 \qquad 0<r<1.
\]
Write $e(f,f)=\int_\D|\nabla f|^2$ and $\rho=|F'|^2$. If $Z$ is reflecting Brownian
motion in $\D$, put
\[
 \Theta_s^{(r)}=\int_0^s\rho_r(Z_t)\,dt,
 \qquad \eta_t^{(r)}=\inf\{s:\Theta_s^{(r)}>t\}.
\]
The process $Z_{\eta_t^{(r)}}$ has Dirichlet energy $e$ and speed
measure $\rho_r\,dx$. Its image under $F_r$ is reflecting Brownian
motion in $G_r$ \cite[Section~6.2]{FOT2011}.
For a Lipschitz domain, $\rho\in L^p(\D)$ for some $p>1$ and
$\rho_r\to\rho$ in $L^p(\D)$. The conformal-exhaustion argument in
\cite[Lemmas~4.1--4.3]{LipschitzCompanion} uses this convergence to obtain spectral
convergence and, when the limiting eigenvalue is simple, convergence
of the corresponding eigenfunction.

\end{document}